\documentclass[a4paper]{article}
\usepackage[frenchb]{babel}
\usepackage[latin1]{inputenc}
\usepackage{amsmath}
\usepackage{array}
\usepackage{amsthm}
\usepackage{amssymb}
\input xy
\xyoption{all}

\author{Aur\'elien DJAMENT\thanks{djament@math.univ-paris13.fr}\;\thanks{http://www.math.univ-paris13.fr/~djament/}\;\thanks{LAGA, Institut Galil\'ee,
universit\'e Paris 13,
99 avenue J.-B. Clément,
93430 VILLETANEUSE, FRANCE}}

\title{Cat\'egories de foncteurs en grassmanniennes et filtration de Krull}

\date{Novembre 2006}

\newcommand{\A}{{\mathcal{A}}}

\newcommand{\C}{{\mathcal{C}}}

\newcommand{\F}{{\mathcal{F}}}
\newcommand{\E}{{\mathcal{E}}}

\newcommand{\K}{{\mathcal{K}}}
\newcommand{\N}{{\mathcal{N}}}
\newcommand{\FF}{{\mathbb{F}_2}}

\newcommand{\Gr}{{\mathcal{G}r}}

\newcommand{\kk}{{\Bbbk}}

\newcommand{\ptt}{{\,\widetilde{\otimes}\,}}

\newtheorem{thm-intro}{Théorème}
\newtheorem{cor-intro}[thm-intro]{Corollaire}
\newtheorem{conj-intro}[thm-intro]{Conjecture}
\newtheorem{pr-intro}[thm-intro]{Proposition}

\newtheorem{theo}{Théorème}[section]
\newtheorem{pr}[theo]{Proposition}
\newtheorem{cor}[theo]{Corollaire}
\newtheorem{lm}[theo]{Lemme}

\newtheorem{prdef}[theo]{Proposition et définition}
\newtheorem{conj}[theo]{Conjecture}

\theoremstyle{definition}
\newtheorem{defi}[theo]{Définition}
\newtheorem{nota}[theo]{Notation}
\newtheorem{hyp}[theo]{Hypothèse}
\newtheorem{conv}[theo]{Convention}

\newtheorem{def-intro}[thm-intro]{Définition}

\theoremstyle{remark}
\newtheorem{rem}[theo]{Remarque}

\begin{document}


\maketitle

\bigskip

\noindent
\textit{Résumé : }   Soit $\F$ la catégorie des foncteurs entre espaces
  vectoriels sur le corps à deux éléments. \`A l'aide des {\em catégories de foncteurs en
  grassmanniennes}, nous avons émis dans \cite{art2} une conjecture
décrivant la filtration de Krull de la catégorie $\F$. Nous démontrons
une forme affaiblie de cette conjecture, avec comme application la
détermination de la structure du produit tensoriel entre le foncteur
projectif standard $P_{\FF^2}$ et un foncteur fini de $\F$, dont on
établit le caractère noethérien. 
Nous étudions également
le morphisme induit par le foncteur d'intégrale en grassmanniennes
entre anneaux de Grothendieck.
\bigskip

\bigskip

\noindent
\textit{Abstract \textbf{(Grassmannian functor categories and Krull filtration)} : } Let $\F$ be the category of functors between
  vector spaces over the field $\FF$. With the help of {\em grassmannian functor
    categories}, we formulated in \cite{art2} a conjecture which describes
  the Krull filtration of the category $\F$. We prove a weak form of
  this conjecture, which applies to show the noetherian character of the tensor product between the standard projective
  functor $P_{\FF^2}$ and a finite functor of $\F$. We study also the
  morphism induced by the grassmannian integral functor between
  Grothendieck rings.

\bigskip

\bigskip

\noindent
\textit{Classification mathématique par sujets : } 16P60, 18A25, 20C33 ; secondaire : 16E20, 16P40, 18E35,
  18G05, 55S10.
\bigskip

\noindent
\textit{Mots clefs : } Catégories de foncteurs, filtration de Krull, objets noethériens, algèbre homologique, foncteur différence et
  filtration polynomiale,  grassmanniennes, représentations modulaires.

\pagebreak

\tableofcontents

\pagebreak

\section*{Introduction}

Soient $\E$ la catégorie des espaces vectoriels sur le corps $\FF$ à
deux éléments, $\E^f$ la sous-catégorie pleine des espaces de
dimension finie et $\F=\mathbf{Fct}(\E^f,\E)$ la catégorie des
foncteurs de $\E^f$ vers $\E$. Cette catégorie joue un rôle important
en algèbre et en topologie ; de nombreuses
investigations sur sa structure ont été menées (cf. \cite{FFPS}). Pour autant, on ignore si $\F$ est une catégorie
localement noethérienne. Ce problème connu sous le nom de {\em
  conjecture artinienne}, discuté dans \cite{GP1} et \cite{art2} par
exemple, n'est qu'une illustration des questions qui
demeurent ouvertes dans le domaine. 

Cet article expose des avancées sur la conjecture artinienne et
précise sensiblement les obstacles à surmonter pour parvenir à une
compréhension globale satisfaisante de la catégorie $\F$. Il combine deux outils essentiels  dans l'étude de la
catégorie $\F$ : les foncteurs $\tilde{\nabla}_n : \F\to\F$ de Powell
et les {\em catégories de foncteurs en grassmanniennes}, introduites dans
\cite{art2}, dont ce travail constitue la continuation. Tous les
termes ou notations non définis ici sont introduits dans la première
partie du présent article, ou
à défaut dans \cite{art2}.

Définis dans \cite{GP4}, les foncteurs $\tilde{\nabla}_n$ constituent une filtration décroissante de
sous-foncteurs du foncteur différence $\Delta : \F\to\F$ obtenue à
partir de la filtration polynomiale de l'injectif standard
$I_\FF$. Ils ont permis à
Powell de donner les premiers renseignements profonds connus sur la
structure {\em globale} de la catégorie $\F$ : son {\em théorème de
simplicité} (établi dans \cite{GP2}) procure des informations sur {\em
tous} les projectifs standard de la catégorie~$\F$. De plus, Powell a
déterminé la structure du projectif $P_{\FF^2}$, dont il a montré le
caractère noethérien de type~$2$, à l'aide du foncteur
$\tilde{\nabla}_2$ (cf. \cite{GP5}). La préservation des
monomorphismes et des épimorphismes par les foncteurs
$\tilde{\nabla}_n$ et le contrôle précis de leur effet sur les
foncteurs simples de $\F$ permettent de fait de mener efficacement des
raisonnements à la fois explicites (en termes d'éléments) et généraux
dans la catégorie~$\F$.

D'un autre côté, les catégories de foncteurs en grassmanniennes
fournissent un cadre puissant pour mener des calculs cohomologiques
dans la catégorie $\F$, à l'aide du foncteur d'intégrale $\omega
:\F_\Gr\to\F$. On rappelle que $\F_\Gr$ est la catégorie des foncteurs
de but $\E$ et de source la catégorie des objets de $\E^f$ munis d'un
sous-espace, et que $\omega$ est donné sur les objets par 
$$\omega(X)(V)=\bigoplus_{W\subset V} X(V,W).$$
Ainsi, le théorème d'annulation cohomologique principal
de \cite{art2} contient les lemmes techniques employés dans \cite{GP5}
pour contrôler des groupes d'extensions apparaissant dans l'étude du
foncteur $P_{\FF^2}$. De surcroît, le foncteur $\omega$ nous a permis,  à l'aide des catégories $\F_{\Gr,n}$,
de donner une description
conjecturale, la {\em conjecture
artinienne extrêmement forte}, de la filtration de Krull $(\K_n(\F))$ de la catégorie
$\F$, définie par
$\K_{-1}(\F)=\{0\}$ et par le fait que, pour $n\geq 0$,
$\K_{n}(\F)/\K_{n-1}(\F)$ est la plus petite sous-catégorie épaisse et
stable par colimites de $\F/\K_{n-1}(\F)$ qui en contient les objets simples. Le c\oe ur du présent
article est constitué par le résultat suivant, où $(\overline{\N
    il}_{\tilde{\nabla}_n})$ est une
filtration de la catégorie $\F$ par des sous-catégories épaisses
définies à partir des foncteurs $\tilde{\nabla}_n$. 

\begin{thm-intro}[Théorème de simplicité généralisé]\label{thint-tsg}
  Pour tout entier $n\geq 0$, le foncteur $\omega_n : \F_{\Gr,n}\to\F$
  induit une équivalence entre la sous-catégorie pleine des objets
  localement finis de $\F_{\Gr,n}$ et une sous-catégorie épaisse de
  $\overline{\N il}_{\tilde{\nabla}_{n+1}}/\overline{\N il}_{\tilde{\nabla}_n}$.
\end{thm-intro}

Le théorème de simplicité de Powell énonce pour sa part que si $X$ est
un objet simple pseudo-constant de $\F_{\Gr,n}$, alors l'image de
$\omega_n(X)$ dans la catégorie  $\overline{\N
  il}_{\tilde{\nabla}_{n+1}}/\overline{\N il}_{\tilde{\nabla}_n}$ est
simple. Nous gagnons donc essentiellement le passage à un objet simple quelconque de
$\F_{\Gr,n}$ et le contrôle des groupes d'extensions dans la catégorie
quotient à partir des groupes d'extensions dans $\F_{\Gr,n}$.

\smallskip

La démonstration du théorème~\ref{thint-tsg} utilise deux types d'ingrédients.
D'une part apparaissent des considérations explicites relatives aux foncteurs
$\tilde{\nabla}_n$, qui procèdent des mêmes idées que le théorème de
simplicité de Powell. D'autre
part, des propriétés d'annulation cohomologique du foncteur $\omega$,
fournies par le théorème principal de \cite{art2} et une variante en
termes des foncteurs $\tilde{\nabla}_n$ établie dans la
section~\ref{sct-nnil}, jouent un rôle significatif.

Une grande part de  la riche structure de la catégorie de
foncteurs en grassmanniennes $\F_\Gr$ intervient dans le théorème de
simplicité généralisé. Ainsi, la filtration par les sous-catégories
$\F_{\Gr,\leq n}$ est omniprésente ; les propriétés cohomologiques du
foncteur $\omega$ utilisées dans la démonstration reposent sur la description fonctorielle de
ces catégories. 

On emploie également de manière décisive la
description en termes de comodules des catégories de foncteurs en grassmanniennes, dans un argument de stabilisation
qui constitue la partie la plus concrète de la démonstration. Cet
argument généralise de manière conceptuelle des considérations déjà
utilisées par Powell. Quant à
la description monadique des catégories $\F_{\Gr,n}$, elle apparaît en filigrane dans les estimations des
foncteurs composés $\tilde{\nabla}_n\,\omega_n$ (pour $n=1$, c'est
essentiellement le début de la résolution canonique qui intervient).

\smallskip

Comme conséquence du théorème de simplicité généralisé et des
résultats de Powell sur les foncteurs annihilés par $\tilde{\nabla}_2$
(utilisés pour déterminer la structure de $P_{\FF^2}$), nous donnons
la contribution suivante à la conjecture artinienne :

\begin{thm-intro}\label{thintc} Pour tout foncteur fini $F$ de $\F$,
  le foncteur $P_{\FF^2}\otimes F$ est noethérien de type~$2$
  (i.e. noethérien et dans $\K_2(\F)$).
\end{thm-intro}

Plus généralement, grâce au théorème~\ref{thint-tsg}, la conjecture
artinienne extrêmement forte est ramenée à un problème plus concret :
montrer que les quotients de la filtration par
$\tilde{\nabla}_n$-nilpotence de certains foncteurs sont oméga-adaptés
d'une certaine hauteur, selon la terminologie introduite dans
\cite{art2}. Cette filtration est définie par le
noyau de flèches entièrement explicites (cf. \cite{GP4}) ; son étude
se trouve étroitement liée à des questions fines de représentations
des groupes symétriques ou linéaires (sur~$\FF$). En effet, les
partitions associées aux facteurs de composition des quotients de
cette filtration sont contrôlées, le problème consiste donc d'une
certaine façon à
maîtriser les extensions entre les différents facteurs de composition possibles. Ce problème dérive essentiellement de la théorie des
représentations ; c'est une sorte de générisation de la question
(beaucoup plus élémentaire) suivante : déterminer les sous-représentations maximales
d'un produit tensoriel de puissances extérieures de la représentation
régulière d'un groupe symétrique dont tous les facteurs de composition
correspondent à des partitions de longueur inférieure à un entier
fixé.

\smallskip

Un autre résultat issu de l'étude conjointe des foncteurs
$\tilde{\nabla}_n$ et $\omega$, donné dans la section~\ref{sctfco}, s'exprime en termes de groupes de
Grothendieck. On peut le formuler explicitement comme suit. Les facteurs de
composition sont implicitement comptés avec multiplicité. On
remarquera qu'un foncteur de type fini de $\F$ (ou plus généralement, prenant
des valeurs de dimension finie) a en général un nombre
infini de facteurs de composition, mais que chacun y a une multiplicité finie.

\begin{thm-intro}\label{thintg} Soient $X$ et $Y$ deux objets finis de la catégorie
  $\F_\Gr$. Supposons que les objets $\omega(X)$ et $\omega(Y)$ de
  $\F$ ont les mêmes facteurs de composition. Alors $X$ et $Y$ ont les mêmes facteurs de composition.
\end{thm-intro}

Ce théorème revêt une double importance. D'une part, il donne une
description du groupe de Grothendieck $G_0^{tf}(\F)$ des objets de type fini de la
catégorie $\F$, en admettant la conjecture artinienne extrêmement
forte. On en déduit ainsi le corollaire ci-dessous, où les catégories de
foncteurs en grassmanniennes n'apparaissent pas. Ce corollaire
illustre la puissance de la conjecture artinienne extrêmement forte,
puisque le résultat qu'elle implique semble de la plus haute
difficulté à atteindre même en admettant les formes renforcées de la
conjecture artinienne  émises avant la forme extrêmement forte.

\begin{cor-intro} Si la conjecture artinienne extrêmement
forte est vraie, alors la classe d'un objet de type fini de $\F$ dans le
groupe de Grothendieck $G_0^{tf}(\F)$ est déterminée par ses facteurs
de composition.
\end{cor-intro}

D'autre part, l'intérêt de la détermination de facteurs de composition
significatifs dans des foncteurs du type $\omega(X)$, où $X$ est un
objet fini de $\F_\Gr$, dépasse le seul cadre des groupes de
Grothendieck : des considérations analogues sont utilisées dans
l'article \cite{art2}, qui établit le théorème~\ref{thintc} dans le
cas particulier des puissances extérieures, sans recours aux
catégories de foncteurs en grassmanniennes. Plus précisément, le
principe du théorème~\ref{thintg}, similaire à certains arguments de
\cite{art2}, consiste à étudier des facteurs de
composition~\guillemotleft~significatifs~\guillemotright~et
raisonnablement détectables dans des
foncteurs du type $\omega(X)$, où $X$ est un objet fini de $\F_\Gr$. Il s'agit de facteurs de composition
associés à des partitions de longueur maximale présentant une
certaine périodicité (une partition s'obtenant à partir d'une
autre en ajoutant ou retranchant un entier à toutes ses parts). En
effet, les foncteurs $\tilde{\nabla}_n$ sont adaptés à l'étude de ces
facteurs de composition d'un foncteur $\omega(X)$, et ceux-ci
caractérisent essentiellement l'objet $X$ de $\F_\Gr$.

Les problèmes de représentations
modulaires auxquels  est liée la démarche du
théorème~\ref{thintg} participent de la richesse de la
structure globale de la catégorie~$\F$. De fait, le résultat
d'injectivité du théorème~\ref{thintg} s'obtient par des arguments
qualitatifs hautement non explicites, de sorte que les facteurs de
composition des foncteurs $\omega(X)$ restent difficiles à étudier de
manière générale, même en se limitant aux
facteurs~\guillemotleft~significatifs~\guillemotright~évoqués plus haut.

\smallskip

Ce travail ne considère que le corps à deux éléments, d'une part pour en
simplifier la partie technique (notamment en évitant le recours à la
décomposition scalaire), d'autre part parce que les considérations non
formelles relatives aux foncteurs annihilés par $\tilde{\nabla}_2$ ne
se généralisent pas clairement : nous ne pouvons pas encore démontrer
l'analogue du théorème~\ref{thintc} sur un corps fini arbitraire. 

La plupart des résultats de cet article sont contenus dans la thèse de
doctorat de l'auteur (\cite{these}).

\subsection*{Notations, conventions et rappels utilisés dans tout
  l'article} {\em Les notations et conventions générales de l'article
  \cite{art2} sont conservées} ; nous ne rappelons que les plus
usitées. D'autres  notations sont rappelées dans la partie~\ref{pun}.
\begin{enumerate}\item Une {\em catégorie de Grothendieck} est une
  catégorie abélienne avec générateurs et colimites exactes. Une telle
  catégorie possède des enveloppes injectives.
\item Soit $\C$  une sous-catégorie localisante, i.e. épaisse et stable par
colimites, d'une catégorie de
Grothendieck $\A$. 
\begin{enumerate}\item La catégorie quotient\,\footnote{Nous renvoyons à \cite{Gab} pour ce qui concerne les catégories
abéliennes quotients.} $\A/\C$ est une catégorie de Grothendieck ;
le foncteur canonique $\A\to\A/\C$ est exact et commute aux
colimites. On le notera quelquefois $\A\twoheadrightarrow\A/\C$ sans
plus de précisions.
\item Ce foncteur canonique possède un adjoint à droite, appelé {\em
    foncteur section}.
\item On dit qu'un objet $X$ de $\A$ est {\em
  $\C$-fermé} si ${\rm Ext}^i_\A(C,X)=0$ pour  $C\in {\rm Ob}\,\C$  et $i\in\{0,1\}$.
\item On dit qu'un objet $X$ de $\A$ est {\em
  $\C$-parfait} si ${\rm Ext}^*_\A(C,X)=0$ pour tout objet $C$ de
$\C$. Le cas échéant, le foncteur canonique
$\pi : \A\to\A/\C$ induit un isomorphisme ${\rm Ext}^*_\A(A,X)\xrightarrow{\simeq}
{\rm Ext}^*_{\A/\C}(\pi(A),\pi(X))$ pour tout objet $A$ de
$\A$. Ce fait sera utilisé abondamment dans cet article, sans plus de précision.
\end{enumerate}
\item\label{it-f} La catégorie $\F$ est une catégorie de
  Grothendieck ; elle est monoïdale symétrique. Cela vaut plus
  généralement pour toute catégorie de foncteurs dont la source est
  essentiellement petite et le but est la catégorie d'espaces
  vectoriels $\E$.
\item Soit $\A$ une catégorie monoïdale symétrique et $X$ un objet de
  $\A$. Lorsqu'ils existent, les adjoints à droite et à gauche à
  l'endofoncteur $-\otimes X$ de $\A$ sont notés respectivement
  $\mathbf{Hom}_\A (X,-)$ et $(-:X)_\A$ ; l'indice $\A$ sera souvent
  omis. Ces adjoints sont appelés respectivement {\em foncteur hom
    interne} et {\em foncteur de division} par $X$.
\item On note $\F^f$, $\F^{tf}$ et $\F^{df}$ respectivement les sous-catégories
  pleines des foncteurs finis (i.e. de longueur finie), de type fini et à
  valeurs de dimension finie de $\F$.
\item Dans la catégorie $\F$ (ou plus généralement une catégorie du
  type indiqué dans~\ref{it-f}), le foncteur $\mathbf{Hom}_\F(F,-)$ est
  toujours défini, et $(-:F)_\F$ l'est dès que $F$ appartient à $\F^{df}$. 
\item Si $G$ est un groupe fini et $M$ un $\FF[G]$-module
  fini\,\footnote{Pour ce qui concerne la terminlogie et les
    propriétés élémentaires de théorie des représentations, nous
    renvoyons le lecteur à \cite{CR}.},
$\mathbf{Hom}(M,-)$ et $(-:M)$ existent et sont naturellement isomorphes
au produit tensoriel par le module contragrédient $M^*$ de~$M$.
\item Si $\A$ est une catégorie de Grothendieck, on note $G_0^f(\A)$
  le groupe de Grothendieck des objets finis de~$\A$.
\item L'espace vectoriel $\FF^{\oplus n}$ est noté $E_n$.
\end{enumerate}



\part{Rappels}\label{pun}

Nous rappelons succinctement quelques propriétés utiles des différentes
catégories de foncteurs qui interviendront dans cet article.

\section{La catégorie $\F$ et les endofoncteurs $\tilde{\nabla}_n$}\label{sct-fpo}

 Soit $V$ un espace vectoriel de dimension finie. On rappelle que l'objet projectif
  (resp. injectif)  standard de $\F$ associé à $V$ est
  noté $P_V$ (resp. $I_V$) ; le foncteur de décalage par $V$ se note
  $\Delta_V$, et $\Delta$ désigne le {\em foncteur différence}
  de~$\F$. On a $\Delta_V\simeq\mathbf{Hom}_\F(P_V,-)\simeq
  (-:I_V)_\F$, et $\Delta\simeq\mathbf{Hom}_\F(\bar{P},-)\simeq
  (-:\bar{I})_\F$, où $\bar{P}$ (resp. $\bar{I}$) désigne la partie
  sans terme constant de $P_\FF$ (resp. $I_\FF$).

\subsection{Foncteurs polynomiaux et filtration polynomiale}\label{par-fpol} On rappelle qu'un foncteur $F$ de $\F$ est dit {\em polynomial} s'il
existe $n\in\mathbb{N}$ tel que $\Delta^n F=0$, et {\em analytique}
s'il est colimite de foncteurs polynomiaux. On note $\F^n$ la
sous-catégorie pleine des foncteurs polynomiaux de degré au plus $n$
de $\F$ (i.e. des foncteurs $F$ tels que $\Delta^{n+1}F=0$), et
$\F_\omega$ la sous-catégorie pleine des foncteurs analytiques
de~$\F$. Un résultat de base de la théorie est que tous les foncteurs
finis de $\F$ sont polynomiaux (cf. \cite{K1}).
Les considérations élémentaires que nous rappelons dans la suite de ce
paragraphe sont détaillées dans \cite{GP4}, par exemple.

Le foncteur d'inclusion $\F^n\hookrightarrow\F$ admet un adjoint à
droite noté $p_n : \F\to\F^n$ ; explicitement, $p_n(F)$ est le plus
grand sous-foncteur polynomial de $F\in {\rm Ob}\,\F$ de degré au
plus~$n$.

La filtration polynomiale du foncteur injectif standard $I_\FF$, noté
simplement $I$ par la suite, est  explicite : si l'on note
\begin{equation}\label{eq-dft} t_n=\sum_{l\in E_n^*}[l]\in\FF[{\rm hom} (E_n,\FF)]\end{equation}
et $t_n^* : I\to I_{E_n}$ le morphisme de $\F$ induit, la suite
suivante est exacte :
\begin{equation}\label{eq-presi} 0\to p_{n-1}(I)\to I\xrightarrow{t_n^*} I_{E_n}.
\end{equation}

Le résultat  suivant joue un rôle fondamental dans le
comportement des foncteurs $\tilde{\nabla}_n$ de Powell et des
variantes que nous introduirons.

\begin{lm}\label{lm-filtr-pol} L'espace vectoriel
$(I/p_{n-1}(I))(E_i)$ est de dimension $0$ si $i<n$ et $1$ si $i=n$.
\end{lm}

Dualement, le foncteur d'inclusion $\F^n\hookrightarrow\F$ admet un adjoint à
gauche noté $q_n : \F\to\F^n$ ; il envoie un foncteur sur son plus
grand quotient polynomial de degré au plus~$n$. Nous noterons
$k_n(F)$, pour $F\in {\rm Ob}\,\F$, le noyau de l'épimorphisme
canonique $F\twoheadrightarrow q_{n-1}(F)$. Le cas qui nous intéresse
ici est dual de la filtration polynomiale de $I$ : pour le projectif
standard $P_\FF$, noté simplement $P$ par la suite, on peut décrire
$k_n(P)$ de la façon suivante.

Si $V$ est un objet de $\E^f$ et $W$ un sous-espace de $V$, notons
$$s_W=\sum_{v\in W} [v]\in P(V).$$
Alors $k_n(P)(V)$ est le sous-espace vectoriel de $P(V)$ engendré par
les éléments $s_W$, où $W$ parcourt les sous-espaces de dimension $n$
de~$V$.

\subsection{Les foncteurs simples de $\F$}\label{par-simples} Une manière classique de
décrire les objets simples de la catégorie $\F$ consiste à employer
les partitions d'un entier, qui apparaissent dans la théorie des représentations
des groupes symétriques\,\footnote{Une autre approche, suivie dans
  \cite{K2}, consiste à s'appuyer sur la théorie des représentations
  des groupes linéaires.} (cf. \cite{James} à ce sujet).

\begin{defi}\label{def-partit}\begin{enumerate}\item Une {\em
      partition} est une suite décroissante d'entiers
    $(\lambda_i)_{i\in\mathbb{N}^*}$ qui stationne en $0$.
\item La {\em longueur} $l(\lambda)$ d'une partition $\lambda$ est le plus grand entier  tel que $\lambda_{l(\lambda)}>0$. Si $\lambda$ est identiquement nulle, on convient que $l(\lambda)=0$. Par la suite, on identifiera une partition $\lambda$ et le $n$-uplet $(\lambda_1,\dots,\lambda_n)$ si $n\geq l(\lambda)$.
\item Une partition $\lambda$ est dite {\em
  $2$-régulière} si
$\lambda_i>\lambda_{i+1}$ pour $1\leq i< l(\lambda)$ ; le corps de
base étant fixé à $\FF$, nous parlerons
par la suite simplement de partition {\em régulière}.

Nous désignerons par $\mathfrak{p}$ l'ensemble des partitions régulières.
\item Le {\em degré} d'une partition $\lambda$ est l'entier positif
  $|\lambda|=\underset{i\in\mathbb{N}^*}{\sum}\lambda_i$. Une partition de $n\in\mathbb{N}$ est par définition une partition de degré $n$.
\item Si $\lambda$ et $\mu$ sont deux partitions de même degré, nous
  noterons $\lambda\leq\mu$ si
$$\forall n\in\mathbb{N}^*\qquad\sum_{i=1}^n\lambda_i\leq\sum_{i=1}^n\mu_i.$$
\end{enumerate}
\end{defi}

 On rappelle que les foncteurs simples de $\F$ peuvent se paramétriser
par les partitions régulières, à partir des représentations des
groupes symétriques. Nous suivons ici les notations de \cite{GP3}, par
exemple ; une présentation détaillée des foncteurs simples est donnée
dans \cite{PS1}, mais cet article indexe les
foncteurs simples de manière différente (par la partition {\em duale}
de celle que nous utilisons).

\begin{nota}\label{not-simples} On désigne par $S_\lambda$ le foncteur simple associé à
  une partition régulière $\lambda$. On note $I_\lambda$ l'enveloppe
  injective du foncteur $S_\lambda$.
\end{nota}

Ainsi, $S_\lambda$ est un foncteur polynomial de degré $|\lambda|$.

Une des propriétés fondamentales des objets simples de la catégorie
$\F$ est la suivante :

\begin{pr}\label{pr-corpsdec} Les corps d'endomorphismes des foncteurs simples
  $S_\lambda$ sont réduits à $\FF$ : on dit que $\FF$ est un {\em corps de décomposition} de la catégorie $\F$.
\end{pr}

Cette proposition est démontrée dans \cite{K2} à partir du résultat
classique (que nous utiliserons également) selon lequel $\FF$ est un
corps de décomposition de la catégorie des $GL_n$-modules (on note
$GL_n$ pour $GL_n(\FF)$). 

Une conséquence formelle de la proposition~\ref{pr-corpsdec} est qu'on
peut calculer la multiplicité, notée $m_\lambda(F)$, d'un foncteur simple $S_\lambda$ dans un
foncteur $F$, supposé à valeurs de dimension finie (pour assurer la
finitude de cette multiplicité), par la formule
\begin{equation*}
m_\lambda(F)=\dim_\FF {\rm hom}_\F(F,I_\lambda).
\end{equation*}

\begin{nota}\label{not-vdash}
\'Etant donnés une partition régulière $\lambda$ et un foncteur $F\in
{\rm Ob}\,\F$, nous abrégerons l'assertion {\em $S_\lambda$ est
  facteur de composition} (i.e. sous-quotient) {\em de $F$} en $\lambda\vdash F$.
\end{nota}

Le théorème suivant résume les propriétés fondamentales des facteurs de
composition du produit tensoriel, noté $\Lambda^\lambda$, des
puissances extérieures $\Lambda^{\lambda_i}$. Il est l'analogue de propriétés
classiques des modules de Specht dans le contexte des groupes
symétriques (cf. \cite{James}) ; une démonstration dans le cadre de
la catégorie $\F$ est donnée dans l'article \cite{K2} de Kuhn (on prendra
garde au fait que cet article emploie des conventions différentes des
nôtres dans l'indexation des foncteurs simples).

\begin{theo}
\label{fclambda}
Les facteurs de composition de $\Lambda^\lambda$, où $\lambda$ est une
partition de longueur $r$ de $n\in\mathbb{N}$, sont :
\begin{itemize}
\item les $S_\mu$, où $\mu$ parcourt les partitions régulières de $n$
 telles que $\mu\geq\lambda$,
\item des simples du type $S_\mu$ avec $|\mu|<n$, $l(\mu)<r$, $\mu_1\geq\lambda_1$ et $\mu_{r-1}\leq\lambda_r$.
\end{itemize}

En outre, $S_\lambda$ est facteur de composition unique de $\Lambda^\lambda$.
\end{theo}

Nous rappelons maintenant un résultat fondamental dont on dispose sur
les produits tensoriels.

\begin{defi}\label{concat-part} Soient $\lambda$ et $\mu$ deux
  partitions de longueurs respectives $n$ et $r$. Nous appellerons
  {\em partition concaténée} de $\lambda$ et $\mu$ la partition
  obtenue en réordonnant la suite d'entiers
  $(\lambda_1,\dots,\lambda_n,\mu_1,\dots,\mu_r)$ ; elle sera
  notée $(\lambda,\mu)$.
\end{defi}

\begin{pr}[Kuhn]\label{prck} Soient $\lambda$ et $\mu$  des partitions
  régulières telles que la partition $(\lambda,\mu)$ soit régulière. Alors
 $m_{(\lambda,\mu)}(S_\lambda\otimes S_\mu)=1$.
\end{pr}

Cette proposition est le théorème $6.17.2$ de \cite{K2}.

La proposition suivante, également démontrée dans \cite{K2}, établit
le lien entre les objets simples de $\F$ et les représentations simples des
groupes linéaires.

\begin{pr}[Kuhn]\label{prk-rep} Soient $\lambda$ une partition
  régulière et $n\in\mathbb{N}$.
\begin{enumerate}\item Si $n=\lambda_1$, alors $S_\lambda (E_n)$ est un
  $\mathbb{F}_2[GL_n]$-module à gauche simple ; nous le noterons
  $R_\lambda$. 
\item Les $R_\mu$ forment
un système complet de représentants des  $\mathbb{F}_2[GL_n]$-modules
  simples lorsque $\mu$ parcourt les éléments de $\mathfrak{p}$ tels
  que $\mu_1=n$.
\item Le $\FF[GL_n]$-module $S_\lambda (E_n)$ est nul si
  $\lambda_1>n$, égal à $R_\lambda$ si $\lambda_1=n$ et isomorphe à
  $R_{(n,\lambda_1,\lambda_1,\dots,\lambda_r)}$ (où $r=l(\lambda)$) si $\lambda_1<n$.
\end{enumerate}
\end{pr}

Nous conserverons, dans la suite, la notation $R_\lambda$ ; nous
désignerons par $m_{R_\lambda}(M)$ la multiplicité de ce module simple
dans un $GL_n$-module $M$.

\subsection{Les foncteurs $\tilde{\nabla}_n : \F\to\F$ de Powell} On prendra
garde que les foncteurs notés $\tilde{\nabla}_n$ sont les {\em {\bf
    duaux}} de ceux introduits dans \cite{GP4} par Powell, à qui sont dus tous les
résultats de ce paragraphe : si l'on note $\tilde{\nabla}_n^{Pow}$ les
foncteurs $\tilde{\nabla}_n$ de \cite{GP4}, nos foncteurs
$\tilde{\nabla}_n$ satisfont des isomorphismes de dualité $D\circ
\tilde{\nabla}_n\simeq\tilde{\nabla}_n^{Pow}\circ D$ et $D\circ
\tilde{\nabla}_n^{Pow}\simeq\tilde{\nabla}_n\circ D$, où $D$ désigne
le foncteur de dualité de $\F$. Il n'en résultera aucune confusion par
la suite, car nous n'utiliserons jamais les foncteurs
$\tilde{\nabla}^{Pow}_n$, adaptés à l'étude de foncteurs analytiques,
tandis que nous traiterons de foncteurs co-analytiques.

\begin{defi}\label{def-nablan} \'Etant donné un entier $n>0$, on note $\tilde{\nabla}_n$
l'endofoncteur de $\F$ noyau de la transformation naturelle
$\Delta_\FF\simeq (-:I)\twoheadrightarrow (-:p_{n-1}(I))$ induite par l'inclusion
$p_{n-1}(I)\hookrightarrow I$.
\end{defi}

\begin{rem} Le foncteur $\tilde{\nabla}_n$ préserve les foncteurs à
valeurs de dimension finie et les foncteurs de type fini.
\end{rem}

\begin{pr}\label{evid-nab}\begin{enumerate}\item Le foncteur
    $\tilde{\nabla}_n$ est additif ; il conserve les injections et les surjections.
\item\label{p4n} Deux foncteurs $\tilde{\nabla}_n$ et $\tilde{\nabla}_m$ (où
  $m\in\mathbb{N}$) commutent toujours (à isomorphisme naturel près) ;
  en particulier, le foncteur $\tilde{\nabla}_n$ commute  au foncteur différence.
\end{enumerate}
\end{pr}

\begin{pr}\label{nab-proj} Il existe dans $\F$ un isomorphisme
  $\tilde{\nabla}_n(P_V)\simeq k_n(P)(V^*)\otimes P_V$ naturel en
  l'espace vectoriel de dimension finie $V$. En particulier,
  $\tilde{\nabla}_n(P_{E_n})\simeq P_{E_n}$ et
  $\tilde{\nabla}_n(P_{E_k})=0$ pour $k<n$.
\end{pr}

Les propositions~\ref{evid-nab} et ~\ref{nab-proj} sont démontrées dans \cite{GP4}.

\begin{nota}\label{unoti}\begin{enumerate}
\item Si $\lambda$ est une partition de longueur $r$, nous noterons $\lambda_{-}$ la partition de longueur $\leq r$ définie par $(\lambda_{-})_{i}=\lambda_{i}-1$ pour $1\leq i\leq r$.
\item Si $\alpha=(\alpha_{1},\dots,\alpha_{k})$ est une suite finie
  d'entiers et $a\in\mathbb{N}$, on notera
  $\alpha_{+a}=(\alpha_{1}+a,\dots,\alpha_{k}+a)$. Cette notation est étendue au cas où $\alpha$ est une partition, en prenant $k=l(\alpha)$ (longueur de la partition~$\alpha$).
\end{enumerate}
\end{nota}

L'un des grands intérêts des foncteurs $\tilde{\nabla}_n$, qui
explique leur efficacité dans de nombreux calculs explicites sur les
foncteurs simples, réside dans le résultat suivant, établi dans le §\,4.3
de \cite{GP4}.

\begin{pr}\label{nab-simp} Soit $\lambda$ une partition régulière. 
\begin{enumerate}\item Le foncteur $\tilde{\nabla}_n(S_\lambda)$ est nul si et
  seulement si $l(\lambda)<n$.
\item Supposons $\lambda$ de
  longueur $n$. On a
  $\tilde{\nabla}_n(S_\lambda)\simeq S_{\lambda_-}$.
\end{enumerate}
\end{pr}

\begin{prdef}\label{nabnil-ep} La sous-catégorie pleine des foncteurs
  $\tilde{\nabla}_n$-nilpotents de $\F$ (i.e. des foncteurs $F$ tels
  qu'il existe $k\in\mathbb{N}$ tel que $(\tilde{\nabla}_n)^k (F)=0$)
  est épaisse. On la note $\N
  il_{\tilde{\nabla}_n}$.

Nous désignerons par $\overline{\N il}_{\tilde{\nabla}_n}$ la plus petite
sous-catégorie localisante de $\F$  contenant $\N il_{\tilde{\nabla}_n}$.
\end{prdef}

C'est le théorème $4.2.3$ de \cite{GP4} (l'hypothèse d'analycité imposée
par Powell aux foncteurs $\tilde{\nabla}_n$-nilpotents n'intervient pas
dans la démonstration).

\smallskip

La conjecture suivante, discutée par Powell dans \cite{GP1}, §\,$3$, est une
variante forte de la conjecture artinienne.

\begin{conj}\label{ca-nabla} Un objet de type fini de $\F$ est noethérien de type $n-1$ si et
  seulement s'il est $\tilde{\nabla}_n$-nilpotent.
\end{conj}

Nous terminons ce paragraphe par un résultat concernant les produits
tensoriels. Il est démontré dans \cite{GP-dim} (théorème $1$), article qui précise
en quoi la notion de $\tilde{\nabla}_n$-nilpotence fournit une~\guillemotleft~bonne~\guillemotright~notion de dimension dans $\F$.

\begin{pr}\label{nab-tens} Le produit tensoriel d'un foncteur
  $\tilde{\nabla}_n$-nilpotent et d'un foncteur
  $\tilde{\nabla}_m$-nilpotent (où $m\in\mathbb{N}$) est
  $\tilde{\nabla}_{n+m-1}$-nilpotent. En particulier, le produit
  tensoriel d'un foncteur $\tilde{\nabla}_n$-nilpotent et d'un foncteur
  fini est $\tilde{\nabla}_n$-nilpotent.
\end{pr}

\section{La catégorie $\F_\Gr$ et le foncteur $\omega$}

Nous rappelons ici les définitions et propriétés des catégories de
foncteurs en grassmanniennes que nous utiliserons le plus fréquemment
par la suite. Nous renvoyons le lecteur à \cite{art2} pour les
détails.

On rappelle qu'une convention générale consiste à omettre les indices
ou exposants désignant une partie de $\mathbb{N}$ lorsque celle-ci est
égale à $\mathbb{N}$ tout entier.

\subsection{Définitions}\label{par-deffgr} Une première famille de
catégories et de foncteurs auxiliaires est donnée comme suit.
\begin{defi}\label{def-rap-fsurj}\begin{enumerate}\item La catégorie
    $\E^f_{surj}$ est la sous-catégorie de $\E^f$ qui a  les
    mêmes objets et dont les morphismes sont les épimorphismes de $\E^f$.
\item La catégorie $\F_{surj}$ est la catégorie de foncteurs $\mathbf{Fct}(\E^f_{surj},\E)$.
\item On note $o : \F\to\F_{surj}$ le foncteur d'oubli,
  i.e. de précomposition par l'inclusion
  $\E^f_{surj}\hookrightarrow\E^f$.
\item Le foncteur $\varpi : \F_{surj}\to\F$ est donné sur
  les objets par 
$$\varpi(X)(V)=\bigoplus_{W\in\Gr(V)}X(W).$$
\item Si $I$ est une partie de $\mathbb{N}$, on note $\E^I_{surj}$ la
  sous-catégorie pleine de $\E^f_{surj}$ des espaces vectoriels dont
  la dimension appartient à $I$, et l'on pose $\F^I_{surj}=\mathbf{Fct}(\E^I_{surj},\E)$.
\end{enumerate}
\end{defi}

Lorsque $I$ est une partie réduite à un élément $n$, la catégorie
$\F^n_{surj}$ est équivalente à celle des $\FF[GL_n]$-modules à gauche,
notée $_{\FF[GL_n]}\mathbf{Mod}$ ;
nous utiliserons souvent cette identification tacitement.

Rappelons que l'on dispose d'un foncteur de prolongement par zéro
$_{\FF[GL_n]}\mathbf{Mod}\simeq\F^n_{surj}\to\F_{surj}$ ; l'image par
le foncteur composé $_{\FF[GL_n]}\mathbf{Mod}\to\F_{surj}\xrightarrow{\varpi}\F$
d'un $GL_n$-module
simple $R_\lambda$, où $\lambda$ est une partition régulière telle que
$\lambda_1=n$ (cf. proposition~\ref{prk-rep}) sera notée $Q_\lambda$ et est
appelée {\em foncteur de Powell} associé à~$\lambda$. Un cas
particulier fondamental est celui d'une partition à un terme $(n)$ ;
$Q_{(n)}$ est noté aussi $\bar{G}(n)$.

Nous introduisons maintenant les catégories de foncteurs en grassmanniennes.

\begin{defi}\label{def-rap-fgr}\begin{enumerate}\item La catégorie
    $\E^f_\Gr$ est la catégorie dont les objets sont les couples
    $(V,W)$, où $V\in {\rm Ob}\,\E^f$ et $W\in\Gr(V)$, et dont les
    morphismes $(V,W)\to (V',W')$ sont les applications linéaires $f :
    V\to V'$ telles que $f(W)=W'$.
\item On note $\F_\Gr$ la catégorie  de foncteurs
  $\mathbf{Fct}(\E^f_\Gr,\E)$.
\item Si $I$ est une partie de $\mathbb{N}$, on note
  $\E^f_{\Gr,I}$ la sous-catégorie pleine de $\E^f_\Gr$ dont
  les objets sont les couples $(V,W)$ pour lesquels $\dim W\in I$, et l'on note $\F_{\Gr,I}$ la catégorie $\mathbf{Fct}(\E^f_{\Gr,I},\E)$.
\end{enumerate}
\end{defi}

Les catégories de type $\F_{\Gr,I}$ que nous utiliserons seront
surtout les $\F_{\Gr,n}$ et les $\F_{\Gr,\leq n}$. On rappelle que
celles-ci peuvent se voir comme des sous-catégories épaisses de
$\F_\Gr$ via le foncteur de prolongement par zéro. En général,
lorsqu'il est défini, le
foncteur de prolongement par zéro $\F_{\Gr,I}\to\F_{\Gr,J}$ (pour
$I\subset J\subset\mathbb{N}$) est noté
$\mathcal{P}_{I,J}$, tandis que le foncteur de restriction
$\F_{\Gr,J}\to\F_{\Gr,I}$ est noté $\mathcal{R}_{J,I}$.

Les objets projectifs et injectifs standard de la catégorie $\F_{\Gr,J}$ sont notés à l'aide des symboles
$P^{\Gr,J}_A$ et $I^{\Gr,J}_A$ respectivement, où $A$ est l'objet de
$\E^f_{\Gr,J}$ auquel ils sont associés. Dans la catégorie
$\F_{surj}$, on emploie des notations du type $P^{surj}_V$ et $I^{surj}_V$.

\begin{defi}\label{nivcon} On appelle {\em niveau}
  d'un objet $X$ de $\F_{\Gr,I}$ l'élément
$${\rm niv} (X)=sup\,\{\dim W\,|\,(V,W)\in {\rm
  Ob}\,\E^f_{\Gr,I}\;\,X(V,W)\neq 0\}$$
de $I\cup\{-\infty,+\infty\}$. On dit que $X$ est de {\em niveau fini}
si ${\rm niv} (X)<+\infty$.

On définit de même le {\em coniveau} de $X$ comme l'élément
$${\rm coniv} (X)=inf\,\{\dim W\,|\,(V,W)\in {\rm Ob}\,\E^f_{\Gr,I}\;\,X(V,W)\neq 0\}$$
de $I\cup\{+\infty\}$.
\end{defi}

Nous donnons à présent une liste de foncteurs exacts fondamentaux entre les
catégories précédentes.
Soit $\F\otimes\F^I_{surj}$ la catégorie
$\mathbf{Fct}(\E^f\times\E^I_{surj},\E)$. On note :

\begin{itemize}
\item[\textbullet] $\iota_I : \F\to\F_{\Gr,I}$ le foncteur donné par
  $\iota_I(F)(V,W)=F(V)$ ;
\item[\textbullet] $\kappa_I : \F\to\F_{\Gr,I}$ le foncteur donné par
  $\kappa_I(F)(V,W)=F(V/W)$ ;
\item[\textbullet] $\rho_I :
    \F^I_{surj}\to\F_{\Gr,I}$ le foncteur donné par
    $\rho_I(A)(V,W)=A(W)$ ;
\item[\textbullet] $\varepsilon_I : \F_{\Gr,I}\to\F^I_{surj}$ le foncteur
  donné par $\varepsilon_I(X)(V)=X(V,V)$ ;
\item[\textbullet] $\xi_I :
    \F\otimes\F^I_{surj}\to\F_{\Gr,I}$ le foncteur donné par
    $\xi_I(F)(V,W)=F(V,W)$ ;
\item[\textbullet]  $\theta_I :
    \F\otimes\F^I_{surj}\to\F_{\Gr,I}$ le foncteur donné par
    $\theta_I(F)(V,W)=F(V/W,W)$ ;
\item[\textbullet]  $\sigma_I : \F_{\Gr,I}\to\F\otimes\F^I_{surj}$ le
  foncteur donné par $\sigma_I(X)(A,B)=X(A\oplus B,B)$ ;
\item[\textbullet]  $\omega :
  \F_\Gr\to\F$ le foncteur d'{\em intégrale en grassmanniennes}, donné  par 
$$\omega(X)(V)=\bigoplus_{W\in\Gr(V)}X(V,W).$$
\end{itemize}

On a ainsi des isomorphismes canoniques $\omega\circ\rho\simeq\varpi$ et
$\varepsilon\circ\iota\simeq o$.

Pour $n\in\mathbb{N}$, on note $\omega_n : \F_{\Gr,n}\to\F$ le
foncteur composé du prolongement par zéro $\F_{\Gr,n}\to\F_\Gr$ et de
$\omega : \F_\Gr\to\F$.

Les {\em foncteurs de décalage} de $\F_{\Gr,I}$  sont les
endofoncteurs $\Delta^{\Gr,I}_V$ de cette catégorie, où $V\in {\rm
  Ob}\,\E^f$, donnés par $\Delta^{\Gr,I}_V(X)(A,B)=X(V\oplus A,B)$.

Le {\em foncteur différence} $\Delta^{\Gr,I}$ de $\F_{\Gr,I}$ est donné par le
scindement naturel $\Delta^{\Gr,I}_\FF\simeq\Delta^{\Gr,I}\oplus id$ ;
son noyau égale l'image essentielle de $\rho_I :
\F^I_{surj}\to\F_{\Gr,I}$, ses objets sont appelés {\em foncteurs
  pseudo-constants}.

On définit comme dans $\F$ la notion de {\em foncteur polynomial} ou
{\em analytique}.

\subsection{Propriétés}\label{par-prfgr} Dans la proposition
fondamentale suivante, on note $\mathbf{Comod}$ les catégories de
comodules ; la structure comultiplicative sur $\FF[\Gr]$ est duale de
celle d'une algèbre de Boole. 

\begin{pr}\label{prfig} Soit $n\in\mathbb{N}$.
\begin{enumerate}\item Le foncteur $\omega$ est adjoint à
  gauche à $\iota$. Il induit une équivalence de catégories
  entre $\F_{\Gr}$ et la sous-catégorie
  $\mathbf{Comod}_{\kk[\Gr]}$ de $\F$.
\item Le foncteur $\omega_n$ induit une équivalence de catégories
  entre $\F_{\Gr,n}$ et la sous-catégorie
  $\mathbf{Comod}_{\bar{G}(n)}^{fid}$ de $\F$ des
  $\bar{G}(n)$-comodules {\em fidèles}, i.e. dont la comultiplication
  est injective.
\item On a un isomorphisme
$$\omega_I(X\otimes\iota_I(F))\simeq\omega_I(X)\otimes F$$
naturel en les objets $X$ de $\F_{\Gr,I}$ et $F$ de $\F$, où
$I=\mathbb{N}$ ou~$n$.
\end{enumerate}
\end{pr}

Les deux propositions qui suivent donnent les résultats fondamentaux
sur les objets finis des catégories de foncteurs en grassmanniennes.

\begin{pr}\label{prf-grof} Les foncteurs finis de la catégorie
  $\F_{\Gr,I}$ sont polynomiaux.
\end{pr}

\begin{pr}\label{prcqd} Soit $I$ une partie non vide de $\mathbb{N}$.
\begin{enumerate}\item \'Etant donné un objet $X$ de $\F_{\Gr,I}$, les
assertions suivantes sont équivalentes :
\begin{enumerate}\item l'objet $X$ de $\F_{\Gr,I}$ est simple ;
\item l'objet $\sigma_I(X)$ de $\F\otimes\F^I_{surj}$ est simple ;
\item il existe un objet simple $F$ de $\F$ et un objet simple $R$ de
  $\F^I_{surj}$ tel que $X$ est isomorphe à $\kappa_I(F)\otimes\rho_I(R)$.
\end{enumerate}\item Les foncteurs exacts $\sigma_I : \F_{\Gr,I}\to\F\otimes\F^I_{surj}$ et
  $\theta_I$ induisent des isomorphismes d'anneaux (sans
  unité si $I$ est infini)
entre $G_0^f(\F_{\Gr,I})$ et $G_0^f(\F\otimes\F^I_{surj})\simeq
G_0^f(\F)\otimes G_0^f(\F^I_{surj})$ réciproques l'un de l'autre. 
\item Le foncteur exact $\xi_I : \F\otimes\F^I_{surj}\to\F_{\Gr,I}$ induit un isomorphisme entre les groupes
  $G_0^f(\F\otimes\F^I_{surj})$ et $G_0^f(\F_{\Gr,I})$.
\end{enumerate}
\end{pr}

Nous donnons maintenant des propriétés relatives aux foncteurs de division.

\begin{pr}\label{divrho} Soit $I$ une partie de $\mathbb{N}$.
\begin{enumerate}\item Le foncteur $\rho_I :
  \F^I_{surj}\to\F_{\Gr,I}$ est adjoint à gauche à $\varepsilon_I$.
\item Il
  existe dans $\F_{\Gr,I}$ un isomorphisme canonique $(\rho_I(A) : X)\simeq\rho_I(A :
  \varepsilon_I(X))$ pour $A\in {\rm Ob}\,\F^I_{surj}$, et $X\in {\rm
    Ob}\,\F_{\Gr,I}$ à valeurs de dimension finie.
\end{enumerate}
\end{pr}

\begin{pr}\label{domeg} Il existe dans $\F$ un isomorphisme $\omega (X : \iota
  (F))\simeq (\omega(X) : F)$ naturel en les objets $X$ de $\F_\Gr$ et
  $F$ de $\F^{df}$. On a un résultat analogue dans $\F_{\Gr,\leq n}$.
\end{pr}

Une partie de la structure élémentaire de la catégorie
$\F_\Gr$ dépend de l'isomorphisme naturel suivant :
\begin{equation}\label{eq-scii}
  \iota(I_V)\simeq\bigoplus_{W\in\Gr(V)}I^\Gr_{(V,W)}\qquad (V\in {\rm Ob}\,\E^f).
\end{equation}

\begin{pr}\label{divif} Il existe un isomorphisme
$$(X : \iota(I_V))(A,B)\simeq\underset{im\,(C\hookrightarrow V\oplus
  A\twoheadrightarrow A)=B}{\bigoplus_{C\in\Gr(V\oplus A)}} X(V\oplus
A,C)$$
naturel en les objets $V$ de $\E^f$, $(A,B)$ de $\E^f_\Gr$ et $X$ de
$\F_\Gr$.

De plus, pour tout $W\in\Gr(V)$, le monomorphisme scindé
naturel $(X:I^\Gr_{(V,W)})\hookrightarrow (X : \iota(I_V))$ induit par
l'épimorphisme scindé $\iota(I_V)\twoheadrightarrow I_{(V,W)}^\Gr$
fourni par l'isomorphisme~{\rm(\ref{eq-scii})} identifie
$(X:I^\Gr_{(V,W)})(A,B)$ au sous-espace 
$$\underset{im\,(C\hookrightarrow V\oplus
  A\twoheadrightarrow W)=W}{\underset{im\,(C\hookrightarrow V\oplus
  A\twoheadrightarrow A)=B}{\bigoplus_{C\in\Gr(V\oplus A)}}} X(V\oplus
A,C)$$
de $(X : \iota(I_V))(A,B)$.
\end{pr}

Cette proposition est sous-tendue par l'isomorphisme naturel suivant :
\begin{equation}\label{eq-pinj}I^\Gr_{(A,B)}\otimes I^\Gr_{(A',B')}\simeq\bigoplus_{C\in Gr(B,B')}
I^\Gr_{(A\oplus A',C)},
\end{equation}
où $(A,B)$ et $(A',B')$ sont des objets de $\F_\Gr$, et l'on a noté
$Gr(B,B')$ l'ensemble des $C\in\Gr(B\oplus B')$ tels que les
applications linéaires $C\hookrightarrow B\oplus
  B'\twoheadrightarrow B$ et $C\hookrightarrow B\oplus
  B'\twoheadrightarrow B'$ soient surjectives.

\begin{lm}\label{lm-nabom} Soient $(V,W)$ et $(A,B)$ deux objets de
  $\E^f_\Gr$. L'unique morphisme non nul $i_{(A,B)} : \FF\to I^\Gr_{(A,B)}$ induit, par tensorisation
  par $I^\Gr_{(V,W)}$, un morphisme $$I^\Gr_{(V,W)}\to
  I^\Gr_{(V,W)}\otimes I^\Gr_{(A,B)}\simeq\bigoplus_{W'\in
    Gr(W,B)} I^\Gr_{(V\oplus A,W')}$$ (cf. isomorphisme~{\rm(\ref{eq-pinj})}) dont les composantes $I^\Gr_{(V,W)}\to I^\Gr_{(V\oplus
    A,W')}$ sont induites par la projection $V\oplus A\twoheadrightarrow V$.
\end{lm}

\begin{proof}  La composante $I^\Gr_{(V,W)}\to I^\Gr_{(V\oplus
    A,W')}$ de $i_{(A,B)}\otimes I^\Gr_{(V,W)}$ s'obtient par application du foncteur
  $E\mapsto \FF^E : \mathbf{Ens}^{op}\to\E$ à la transformation
  naturelle
$${\rm hom}_{\E^f_\Gr} (\cdot,(V\oplus A,W'))\hookrightarrow {\rm
  hom}_{\E^f_\Gr} (\cdot,(V,W))\times {\rm hom}_{\E^f_\Gr} (\cdot,(A,B))\twoheadrightarrow {\rm
  hom}_{\E^f_\Gr} (\cdot,(V,W))\,,$$
qui est induite  par la projection $V\oplus A\twoheadrightarrow V$,
d'où le lemme.
\end{proof}

Nous utiliserons le théorème d'annulation cohomologique fondamental de
\cite{art2} via sa conséquence suivante :

\begin{pr}\label{fondcr} Soient $k$ et $n$ deux entiers naturels, $X$
  un objet analytique de $\F_{\Gr,k}$
  et $Y$ un objet quelconque de $\F_{\Gr,n}$.
\begin{enumerate}\item Si $k<n$, alors ${\rm
    Ext}^*_{\mathcal{F}}(\omega_k (X),\omega_n(Y))=0$.
\item Si $k=n$, alors le morphisme  naturel ${\rm
    Ext}^*_{\Gr,n}(X,Y)\to {\rm
    Ext}^*_{\mathcal{F}}(\omega_n (X),\omega_n(Y))$ induit par $\omega_n$ est un isomorphisme.
\end{enumerate}
\end{pr}

\subsection{Foncteurs oméga-adaptés}\label{par-omad} Soit
$n\in\mathbb{N}$. La sous-catégorie
pleine des foncteurs {\em oméga-adaptés de hauteur au plus $n$}  de $\F$, notée $\F^{\omega - ad(n)}$ est la plus petite sous-catégorie
qui contient tous les foncteurs du type $\omega(X)$, où $X\in {\rm
  Ob}\,\F_\Gr$ est fini et de niveau au plus~$n$, et qui vérifie la
propriété suivante : si $0\to F\to G\to H\to 0$ est une suite exacte
de $\F$ telle que deux des foncteurs $F$, $G$ et $H$ appartiennent à
$\F^{\omega - ad(n)}$, il en est de même pour le troisième.

\begin{conj}[Conjecture artininenne extrêmement forte]\label{conj-caef} Pour tout $n\in\mathbb{N}$, la
  sous-catégorie $\F^{\omega - ad(n)}$ de $\F$ est stable par quotients.
\end{conj}

Cet énoncé revient à dire que le foncteur $\omega_n$
induit une équivalence entre la sous-catégorie pleine $\F_{\Gr,n}^{lf}$ des
objets localement finis de $\F_{\Gr,n}$ et le sous-quotient
$\K_n(\F)/\K_{n-1}(\F)$ de la filtration de Krull de~$\F$.

La proposition suivante constitue une variante des arguments de la
section~$12$ de \cite{art2}.

\begin{pr}\label{rec-omad} Soit $n\in\mathbb{N}$ tel que les
  sous-catégories $\F^{\omega - ad(i)}$ de $\F$ sont stables par
  quotient pour $i\leq n$. Alors le foncteur $\omega_n$ induit une
  équivalence entre la sous-catégorie pleine $\F^f_{\Gr,n}$ des objets
  finis de $\F_{\Gr,n}$ et la catégorie quotient $\F^{\omega -
    ad(n)}/\F^{\omega - ad(n-1)}$. De plus, si $X$ est un objet de
  $\F^f_{\Gr,n}$, alors $\omega_n(X)$ est un objet noethérien de
  type~$n$ de~$\F$.
\end{pr} 

En termes des foncteurs $\tilde{\nabla}_n$, la conjecture naturelle
relative aux foncteurs oméga-adaptés est la suivante.

\begin{conj}\label{canabla} Un objet de type fini de $\F$ est oméga-adapté de
  hauteur au plus $n-1$ si et
  seulement s'il est  $\tilde{\nabla}_n$-nilpotent.
\end{conj}

Cette conjecture est équivalente à la conjecture artinienne extrêmement forte. Cela
provient de l'épaisseur de $\N il_{\tilde{\nabla}_n}$, de la non
$\tilde{\nabla}_n$-nilpotence des foncteurs $\omega_n(X)$ pour
$X\in {\rm Ob}\,\F_{\Gr,n}$ non nul et de la
$\tilde{\nabla}_n$-nilpotence des foncteurs $\omega_i(X)$ pour $i<n$
et $X\in
{\rm Ob}\,\F_{\Gr,i}$ fini, que nous démontrerons à la
section~\ref{par-nablom}. 

Les conjectures~\ref{conj-caef} et~\ref{canabla} impliquent également la conjecture~\ref{ca-nabla}.

\part{La catégorie $\F/\F_\omega$}\label{pdx}

Cette partie démontre une forme partielle de la conjecture artinienne
extrêmement forte, pour le quotient $\K_1(\F)/\K_0(\F)$
(section~\ref{sct-fom}), qui est ensuite discutée dans la
section~\ref{sct-rqconj}.

L'essentiel des résultats de ce paragraphe ont été déjà établis par un
biais différent par
Powell, dans \cite{GP3}, généralisant les
théorèmes obtenus par Piriou dans \cite{Piriou} à l'aide de méthodes
plus directes.

La démarche présentée dans cette partie présente un double intérêt :
d'une part, elle précise les résultats de Powell,
d'autre part, elle constitue une introduction à la démonstration du
théorème de simplicité généralisée (théorème~\ref{thsg}), dont elle contient
toutes les idées conceptuelles, hormis l'emploi des foncteurs
$\tilde{\nabla}_n$, remplacés par le foncteur exact $\Delta$, ce qui
en simplifie nettement la partie technique.

\section{Le théorème principal}\label{sct-fom}

On rappelle que $\F^{lf}_{\Gr,1}=\K_0(\F_{\Gr,1})$ est la sous-catégorie pleine des objets
localement finis de $\F_{\Gr,1}$. Le fait que les
  objets finis de $\F_{\Gr,1}$ soient de présentation finie
  (cf. \cite{art2}) implique en effet
  formellement que la
  sous-catégorie des objets localement finis de $\F_{\Gr,1}$ est
  épaisse ; de même, $\F_\omega=\K_0(\F)$.

Le but de
cette section est d'établir le résultat suivant :

\begin{theo}\label{th-om1} Le foncteur 
$$\overline{\omega}_1 : \F^{lf}_{\Gr,1}\hookrightarrow\F_{\Gr,1}\xrightarrow{\omega_1}\mathcal{F}\twoheadrightarrow\mathcal{F}/\mathcal{F}_\omega$$
induit une équivalence entre la catégorie $\F^{lf}_{\Gr,1}$ et une sous-catégorie localisante de $\K_1(\F)/\K_0(\F)=\K_1(\F)/\F_\omega$. En particulier, il envoie un foncteur simple de $\mathcal{F}_{\mathcal{G}r,1}$ sur un objet simple de $\mathcal{F}/\mathcal{F}_\omega$.
\end{theo}

En termes de foncteurs oméga-adaptés, ce théorème prend la forme
suivante :

\begin{cor}\label{th1-oad} La sous-catégorie $\F^{\omega - ad(1)}$ de
  $\F$ est épaisse ; ses objets sont des foncteurs noethériens de
  type~$1$.

De plus, le foncteur $\omega_1$ induit une équivalence entre la
sous-catégorie pleine des objets finis de $\F_{\Gr,1}$ et la catégorie
$\F^{\omega - ad(1)}/\F^f$.
\end{cor}

Le théorème~\ref{th-om1} redonne les résultats principaux de l'article
\cite{GP3} de Powell et les précise. Nous y reviendrons dans la section~\ref{sct-rqconj}.

\begin{rem}\label{ca1rq} Il semble en revanche très difficile de
  prouver que tous les objets simples de
  $\mathcal{F}/\mathcal{F}_\omega$ sont oméga-adaptés  de hauteur~$1$ sans avoir démontré
  une version forte de la conjecture artinienne.
\end{rem}

\smallskip

Avant de démontrer le théorème \ref{th-om1}, nous établissons des
résultats préliminaires décrivant le comportement du foncteur
différence sur l'image par $\omega_1$ d'un foncteur fini de
$\F_{\Gr,1}$.

\begin{nota} Soit $X$ un objet de $\F_{\Gr,1}$. Nous désignerons par
  $\pi_X : \Delta\omega_1(X)\to\omega_1(X)$ le morphisme donné sur $V\in {\rm
    Ob}\,\E^f$ par la composée
$$\Delta\omega_1(X)(V)\hookrightarrow\widetilde{\Delta}\omega_1(X)(V)=\omega_1(X)(V\oplus\mathbb{F}_2)=\bigoplus_{(v,t)\in (V\oplus\mathbb{F}_2)\setminus\{0\}}X(V\oplus\mathbb{F}_2,(v,t))$$
$$\twoheadrightarrow\bigoplus_{v\in V\setminus\{0\}}X(V\oplus\mathbb{F}_2,(v,1))\twoheadrightarrow\bigoplus_{v\in V\setminus\{0\}}X(V,v)=\omega_1(X)(V),$$
où la dernière flèche est induite par la projection
$V\oplus\FF\twoheadrightarrow V$.

Si $k$ est un entier naturel et $X$ un
  objet de $\mathcal{F}_{\mathcal{G}r,1}$, nous noterons $\pi_{X}^k
  : \Delta^k\omega_1(X)\to\omega_1(X)$ le morphisme
  $\pi_{X}\circ\Delta\pi_{X}\circ\dots\circ \Delta^{k-1}\pi_{X}$ ; par
  convention, $\pi^0_X=id_X$.
\end{nota}

Ces notations seront utilisées uniquement dans cette section.

Dans la suite de cette section, la catégorie $\F\otimes\F^1_{surj}$ but de $\sigma_1$
est identifiée à~$\F$.

\begin{lm}\label{lm-elm} Soit $X$ un objet de $\mathcal{F}_{\mathcal{G}r,1}$.
\begin{enumerate}
\item\label{lp1} On a un isomorphisme naturel
  $\Delta\omega_1(X)\simeq\sigma_1(X)\oplus\omega_1(\Delta^{\Gr,1} X)\oplus\omega_1\iota_1\sigma_1(X)$.
\item\label{lp2} \`A travers cet isomorphisme, $\pi_X : \Delta\omega_1(X)\to\omega_1(X)$
  se lit comme la composée de la projection
  $\sigma_1(X)\oplus\omega_1(\Delta^{\Gr,1}
  X)\oplus\omega_1\iota_1\sigma_1(X)\twoheadrightarrow\omega_1\iota_1\sigma_1(X)$ et du morphisme obtenu en appliquant $\omega_1$ à la coünité de l'adjonction $\iota_1\sigma_1(X)\to X$.
\item\label{lp3} Le morphisme $\pi_X$ est surjectif.
\item\label{lp4} Si $X$ est fini de degré $d\geq 0$, $ker\,\pi_X$ est la
  somme directe d'un objet fini de $\F$ et de l'image par $\omega_1$
  d'un objet fini  de $\F_{\Gr,1}$ de degré strictement inférieur à~$d$. 
\item\label{lp5} Soient $V$ un espace vectoriel de dimension finie,
  $l$ une forme linéaire non nulle sur $V$ et $x=(x_v)_{v\in
    V\setminus\{0\}}$ ($x_v\in X(V,v)$) un élément de
  $\omega_1(X)(V)$. Notons $a_x : P_V\to\omega_1(X)$ le morphisme de
  $\mathcal{F}$ représenté par $x$. Alors l'élément $(y_v)$ de $\omega_1(X)(V)$ représenté par le morphisme
$$P_V\xrightarrow{u\mapsto [l]\otimes u}\mathbb{F}_2[V^*\setminus\{0\}]\otimes P_V\simeq\Delta P_V\xrightarrow{\Delta a_x}\Delta\omega_1(X)\xrightarrow{\pi_X}\omega_1(X)$$
est donné par $y_v=l(v)x_v$.
\end{enumerate}
\end{lm}

\begin{proof} La
  décomposition ensembliste $$\Gr_1(V\oplus\FF)=\{(0,1)\}\sqcup\{(v,0)\,|\,v\in
  V\setminus\{0\}\}\sqcup\{(v,1)\,|\,v\in
  V\setminus\{0\}\}$$ et l'isomorphisme naturel
$$\bigoplus_{v\in V\setminus\{0\}}X(V\oplus\FF,(v,1))\simeq
\bigoplus_{v\in V\setminus\{0\}}X(V\oplus\FF,(0,1))$$
induit par les isomorphismes
$(V\oplus\FF,(v,1))\xrightarrow{\simeq}(V\oplus\FF,(0,1)),\quad
(x,t)\mapsto (x+tv,t)$ fournissent  un isomorphisme naturel
$$\widetilde{\Delta}\omega_1(X)\simeq\sigma_1(X)\oplus\omega_1(\widetilde{\Delta}^{\Gr,1}
  X)\oplus\omega_1\iota_1\sigma_1(X).$$ L'inclusion canonique
  $\omega_1(X)\hookrightarrow\widetilde{\Delta}\omega_1(X)$ se lit
  dans cet isomorphisme comme $\omega_1(X)\hookrightarrow\omega_1(\widetilde{\Delta}^{\Gr,1}
  X)\hookrightarrow\sigma_1(X)\oplus\omega_1(\widetilde{\Delta}^{\Gr,1}
  X)\oplus\omega_1\iota_1\sigma_1(X)$, d'où le premier point.
 
Pour tout $v\in V\setminus\{0\}$, la composée
$X(V\oplus\FF,(0,1))\xrightarrow{\simeq}X(V\oplus\FF,(v,1))\to
X(V,v)$, où la première flèche est l'inverse de l'isomorphisme utilisé
précédemment et la seconde est induite par la projection, est induite
par $V\oplus\FF\to V,\quad (x,t)\mapsto x+tv$. Cela montre 
l'assertion~\ref{lp2}.

La suite exacte
$\iota_1\sigma_1\Delta^{\Gr,1}(X)\to\iota_1\sigma_1(X)\to X\to 0$ (début de la résolution canonique de $X$ ---
cf. \cite{art2}, §\,$7.2$), et l'exactitude du foncteur $\omega_1$
permettent d'en déduire les assertions~\ref{lp3} et~\ref{lp4}, puisque
les foncteurs $\iota_1$ et $\sigma_1$ respectent le degré.

Pour établir la dernière assertion, notons $z$ l'élément de $\Delta\omega_1(X)(V)\subset\omega_1(X)(V\oplus\mathbb{F}_2)$ représenté par le morphisme
$$P_V\xrightarrow{u\mapsto [l]\otimes u}\mathbb{F}_2[V^*\setminus\{0\}]\otimes P_V\simeq\Delta P_V\xrightarrow{\Delta a_x}\Delta\omega_1(X).$$
On a $z=\big(\omega_1 (X)(i)\big)(x)+\big(\omega_1 (X)(i_l)\big)(x)$, où $i : V\to V\oplus\mathbb{F}_2$ est l'inclusion canonique et $i_l : V\to V\oplus\mathbb{F}_2$ le morphisme de composantes $id_V$ et $l$.

On en déduit que la composante $z_{v,1}$ ($v\in V\setminus\{0\}$) de $z$ égale $X(i_l)(z_v)$ si $l(v)=1$, $0$ sinon. La conclusion provient donc de ce que la composition $V\xrightarrow{i_l}V\oplus\mathbb{F}_2\twoheadrightarrow V$ est l'identité.
\end{proof}

La proposition suivante fournit un argument de stabilisation décisif
relatif aux morphismes $\pi^k_X$.

\begin{pr}\label{cme1} Soient $X$ un objet de $\F_{\Gr,1}$ et $F$ un
  sous-objet de $\omega_1(X)$. Pour tout entier $k\geq 0$, on note
  $C_k=\pi_X^k(\Delta^{k}F)$.
\begin{enumerate}\item La suite $(C_k)_{k>0}$ de sous-objets de $\omega_1(X)$
  est croissante ; nous noterons $C_\infty$ sa réunion. Si $X$ est un objet
  noethérien de $\F_{\Gr,1}$, cette suite stationne. 
\item Pour tout $k\geq 0$, le foncteur $C_k$ est engendré par les
  éléments du type
$$(l_1(v)\dots l_k(v)x_v)_{v\in V\setminus\{0\}}\in\omega_1(X)(V),$$
 où $V$ parcourt les espaces vectoriels de
  dimension finie, $x=(x_v)_{v\in V\setminus\{0\}}$ les éléments de
  $F(V)$ et $(l_1,\dots,l_k)$ les $k$-uplets de formes linéaires sur $V$.
\item Le foncteur $C_\infty$ est le plus petit sous-$\bar{P}$-comodule de
  $\omega_1(X)$ contenant $F$. Si $F$ est lui-même un sous-$\bar{P}$-comodule de
  $\omega_1(X)$, on a $C_k=F$ pour tout $k\geq 0$.
\item Si $X$ est localement fini, on a ${\rm hom}\,(C_\infty/F,\omega_1(U))=0$ pour tout objet $U$ de~$\F_{\Gr,1}$.
\end{enumerate}
\end{pr}

\begin{proof}L'assertion~\ref{lp5} du lemme~\ref{lm-elm} et la préservation
des épimorphismes par le foncteur différence montrent le second point
pour $k=1$. Le cas général s'en déduit aussitôt par récurrence.

On en déduit que la suite $(C_k)_{k>0}$ est croissante. Si $v$ est un
élément non nul d'un espace vectoriel de dimension finie $V$, soient
$l_1,\dots,l_n$ les éléments de $V^*$ tels que $l_i(v)=1$. Alors la
fonction polynomiale $l_1\dots l_n : V\to\FF$ est égale à
l'indicatrice de $\{v\}$. Si l'on note $F^{gr}(V)$ le plus petit sous-espace vectoriel $V\setminus\{0\}$-gradué de $\omega_1(X)(V)$
contenant $F(V)$, ce qui précède montre que $C_\infty$ est le plus petit
sous-foncteur de $\omega_1(X)$ tel que $F^{gr}(V)\subset C_\infty(V)$
pour tout $V\in {\rm Ob}\,\E^f$ ; en particulier, $C_\infty\supset F$. 

Soit $Y$ le
plus petit sous-objet de $X$ tel que $Y(V,v)$ contienne les composantes
dans $X(V,v)$ des éléments de $F(V)\subset\omega_1(X)(V)$. La
proposition~\ref{prfig} montre d'une part que $\omega_1(Y)$ est le plus petit sous-$\bar{P}$-comodule de
  $\omega_1(X)$ contenant $F$. Le paragraphe précédent de la
  démonstration montre d'autre part que
  $C_\infty=\omega_1(Y)$, d'où le troisième point.

Soit $f : C_\infty/F\to\omega_1(U)$ un morphisme de $\F$. Si $X$ est
localement fini, il en est de même pour $Y$, donc la proposition~\ref{fondcr} montre que la composée $g : \omega_1(Y)=C_\infty\twoheadrightarrow
C_\infty/F\xrightarrow{f}\omega_1(U)$ est induite par un morphisme $u : Y\to U$
de $\F_{\Gr,1}$. Comme la composée
$F\hookrightarrow\omega_1(Y)\xrightarrow{g}\omega_1(U)$ est nulle et
que $\omega_1$ est exact, $F$ est inclus dans
$\omega_1(ker\,u)$. D'après le troisième point, on en déduit $ker\,u=Y$,
donc $u=0$, puis $g=0$ et $f=0$, d'où la dernière assertion.

Par ailleurs, si $X$ est noethérien, alors $Y$ est de type fini, donc $C_\infty=\omega_1(Y)$
est de type fini. Cela montre que la suite $(C_k)_{k>0}$ est
stationnaire, et achève la démonstration.
\end{proof}

\begin{proof}[Début de la démonstration du théorème~\ref{th-om1}] La proposition~\ref{fondcr} entraîne la pleine fidélité de $\overline{\omega}_1 : \F^{lf}_{\Gr,1}\to\F/\F_\omega$ et
la stabilité par extensions de son image. Comme le foncteur
$\overline{\omega}_1$ commute aux colimites, il suffit donc
d'établir qu'un sous-objet d'un objet de l'image $\C$ de sa
restriction à la sous-catégorie $\F^f_{\Gr,1}$ des objets finis de
$\F_{\Gr,1}$ est isomorphe à un objet
de $\C$. Pour cela, on procède par récurrence sur le degré polynomial
(cf. proposition~\ref{prf-grof}).

On se donne donc un objet fini $X$
de $\F_{\Gr,1}$ de degré $d\geq 0$ et l'on suppose que l'hypothèse
suivante est satisfaite.

\begin{hyp}[Hypothèse de récurrence] L'image par le foncteur
  $\overline{\omega}_1 : \F^{lf}_{\Gr,1}\to\F/\F_\omega$ du théorème~\ref{th-om1} de la sous-catégorie $(\F^f_{\Gr,1})^{d-1}$
  des foncteurs de degré $<d$ de
  $\F^{f}_{\Gr,1}$ est une sous-catégorie épaisse de
  $\F/\F_\omega$.
\end{hyp}

\begin{nota} Dans cette section, nous noterons $\A_d$ cette sous-catégorie épaisse de
  $\F/\F_\omega$, et $\mathcal{Q}_d$ la catégorie quotient de
$\F/\F_\omega$ par $\A_d$. Nous noterons également $\mathcal{X}$ la sous-catégorie
  pleine des
  objets $T$ de $\F/\F_\omega$ tels que ${\rm
    hom}\,(T,\overline{\omega}_1(U))=0$ pour tout objet $U$ de $\F_{\Gr,1}$.
\end{nota}

\begin{lm}\label{lm-prelf4}\begin{enumerate}\item Le foncteur différence $\Delta$
    induit un endofoncteur exact et fidèle  de la  catégorie
    $\F/\F_\omega$. Il induit également un endofoncteur exact des catégories $\A_d$ et $\mathcal{Q}_d$.
\item La sous-catégorie $\mathcal{X}$ de $\F/\F_\omega$ est stable par le foncteur induit par le
  foncteur différence.  
\item L'intersection des sous-catégories $\mathcal{X}$ et $\A_d$ de
  $\F/\F_\omega$ est réduite à $0$.
\item Pour tout entier $k>0$, le morphisme $\pi^k_X$ induit un
  isomorphisme dans la catégorie $\mathcal{Q}_d$.
\end{enumerate}
\end{lm}

\begin{proof}[Démonstration du lemme] Le foncteur $\Delta$ est exact
  et conserve $\F_\omega$, il induit donc un endofoncteur exact de
  $\F/\F_\omega$. La suite exacte $0\to p_0(F)=F(0)\to F\to\bar{I}\otimes\Delta F$
naturelle en l'objet $F$ de $\F$ montre la fidélité de ce foncteur exact. En effet, le foncteur constant $F(0)$ est analytique, et, comme
$\bar{I}$ est analytique, si $\Delta F$ est objet de $\F_\omega$, il
en est de même pour $\bar{I}\otimes\Delta F$.

Le lemme~\ref{lm-elm} montre que $\Delta$ préserve l'image $\A_d$ de
$(\F^f_{\Gr,1})^{d-1}$ par $\overline{\omega}_1$, il induit donc un
endofoncteur exact de $\A_d$ et
$\mathcal{Q}_d=(\F/\F_\omega)/\A_d$. Cela achève de prouver le premier
point. Le lemme~\ref{lm-elm} montre également que le noyau de l'épimorphisme $\pi_X$ appartient à
$\A_d$, puisque $X$ est supposé de degré $d$. Cela établit,
par récurrence sur $k$, la dernière assertion.

L'isomorphisme d'adjonction ${\rm hom}_\F
(\Delta F,\omega_1(U))\simeq {\rm hom}_\F
(F,\omega_1(U\otimes\iota_1(\bar{I})))$ donne le second point, puisque
la proposition~\ref{fondcr} permet de remplacer les ensembles de morphismes
considérés dans $\F$ par des ensembles de morphismes analogues dans
$\F/\F_\omega$.

Le troisième point résulte de la définition de la sous-catégorie $\A_d$.
\end{proof}

\noindent
{\em Fin de la démonstration du théorème~\ref{th-om1}.} --- Soit $F$ un sous-objet de $\omega_1(X)$ ; on conserve les notations
de la proposition~\ref{cme1}, et l'on se donne $k\in\mathbb{N}^*$ tel
que $C_\infty=C_k$ (qui est donc de la forme $\omega_1(Y)$ pour un sous-objet
$Y$ de $X$). Alors $\Delta^{k}F$ et $\Delta^{k}C_\infty$ ont la même
image $C_\infty$ par $\pi_X^k$, qui induit un isomorphisme dans
$\mathcal{Q}_d$ (dernière assertion du lemme~\ref{lm-prelf4}), donc l'inclusion
$\Delta^{k}F\hookrightarrow\Delta^{k}C_\infty$ induit un isomorphisme
dans $\mathcal{Q}_d$. Ainsi, l'image dans $\F/\F_\omega$ de $\Delta^{k}(C_\infty/F)$ est objet de
$\A_d$. Mais c'est aussi un objet de $\mathcal{X}$ par la dernière assertion de
la proposition~\ref{cme1} et la deuxième assertion du lemme~\ref{lm-prelf4}. La troisième assertion de ce lemme montre alors que
l'image de $\Delta^{k}(C_\infty/F)$ dans $\F/\F_\omega$ est nulle, donc
aussi l'image de $C_\infty/F$ (par la première assertion du lemme). Cela achève la démonstration.
\end{proof}
\section{Remarques et  conjecture}\label{sct-rqconj}
\begin{itemize}
\item[\textbullet] Comme l'article \cite{GP3}, le
  théorème~\ref{th-om1} signifie que les foncteurs $P\otimes F$ (avec $F$ fini) sont
  noethériens de type $1$ et donne des renseignements sur leur
  structure. La méthode qu'emploie Powell dans \cite{GP3} repose
  fortement sur les propriétés des foncteurs $\tilde{\nabla}_n$ et
  certaines considérations explicites sur les représentations des
  groupes symétriques. La nôtre reste très générale et clarifie les
  calculs d'algèbre homologique entre les différents objets dont on
  démontre le caractère simple noethérien de type $1$, entièrement
  ramenés à des calculs d'algèbre homologique sur des objets finis de
  $\F_{\Gr,1}$, qui peuvent théoriquement se comprendre à partir des
  représentations d'algèbres de dimension finie sur~$\FF$. 

Le théorème~\ref{th-om1} fournit également une construction  de l'objet simple noethérien de type $1$ et $\F_\omega$-parfait
associé à une partition régulière $\lambda$ différente de celle de Powell. Nous pouvons aussi montrer, sans calcul, l'égalité
$X_\lambda=K_\lambda$ conjecturée dans \cite{GP3}, §\,$4.5$ : avec nos
notations, elle se réduit à dire que l'image de l'unique morphisme non
nul $\iota_1(S_\lambda)\to\iota_1(S_{\hat{\lambda}})$,
où $\hat{\lambda}=(\lambda_1+1,\dots,\lambda_r+1)$ si $\lambda$ est de
longueur $r$, est $\kappa_1(S_\lambda)$. En effet, cette image
contient $\kappa_1(S_\lambda)={\rm cosoc}\,\iota_1(S_\lambda)$ (cf. \cite{art2}), et
lorsqu'on applique le foncteur exact et fidèle $\sigma_1$ au morphisme
précédent, on obtient l'unique morphisme non nul (cf. théorème B.3 de \cite{PS1})
$S_\lambda\oplus\Delta S_\lambda\to S_{\hat{\lambda}}\oplus\Delta
S_{\hat{\lambda}}$, dont l'image est
$S_\lambda=\sigma_1\kappa_1(S_\lambda)$.
\item[\textbullet] On peut généraliser le théorème~\ref{th-om1} au cas d'un corps fini
quelconque $\kk$. Il convient pour cela de remplacer le foncteur différence
par la composée
$\F(\kk)\xrightarrow{\Delta}\F(\kk)\xrightarrow{(-)_i}\F(\kk)$, où la
dernière flèche provient de la décomposition scalaire, pour
tout entier $i\in\{1,\dots,q-1\}$, où $q={\rm Card}\,\kk$.
\item[\textbullet] Les objets finis de la catégorie $\F/\F_\omega$ donnés par
  le théorème~\ref{th-om1} (il n'y en a pas d'autres si la conjecture
  artinienne extrêmement forte est vraie) possèdent un représentant
  $\F_\omega$-parfait dans $\F$, du type $\omega_1(X)$ (avec $X\in
  {\rm Ob}\,\F_{\Gr,1}$). Il n'en est pas de même pour tous les
  objets de $\F/\F_\omega$. Considérons par exemple le noyau $A$ de
  l'unique morphisme non nul $\bar{P}^{\otimes 2}\to\bar{P}$. \`A partir
  de la suite exacte
$$0\to A\to\bar{P}^{\otimes 2}\to\bar{P}\to\Lambda^1\to 0$$
et de la proposition~\ref{fondcr}, on obtient un isomorphisme naturel
${\rm Ext}^i(F,A)\simeq {\rm Ext}^{i-2}(F,\Lambda^1)$ pour $F\in {\rm
  Ob}\,\F_\omega$. Ainsi, $A$ est $\F_\omega$-fermé, mais non
$\F_\omega$-parfait.
\item[\textbullet] L'absence de représentants parfaits pour certains foncteurs
  équivaut à l'inexactitude du foncteur section
  $\F/\F_\omega\to\F$. Notons néanmoins  la conséquence suivante de
  la proposition~\ref{fondcr}, dans laquelle nous
  notons $r^i$ le $i$-ème foncteur dérivé droit de l'endofoncteur de
  $\F$ composé du
  foncteur canonique $\F\to\F/\F_\omega$ et du foncteur section :
  pour tout foncteur oméga-adapté $F$, on a $r^i(F)=0$ pour $i$ assez
  grand. Ainsi, la conjecture artinienne extrêmement forte permettrait
  de contrôler le défaut d'exactitude du foncteur section.
\item[\textbullet] Même en admettant la conjecture artinienne extrêmement forte,
  certains aspects de la structure globale de la catégorie
  $\F/\F_\omega$ demeurent mystérieux. Par exemple, la description de
  l'enveloppe injective de l'image du foncteur $\bar{P}$ dans cette
  catégorie quotient, qui en est l'objet simple le plus élémentaire,
  pose des problèmes très difficiles.
\end{itemize}

En général, les injectifs indécomposables de $\F/\F_\omega$ sont les
images par le foncteur section des injectifs indécomposables sans
sous-objet fini non nul de $\F$ (cf. \cite{Gab}),
i.e. des injectifs indécomposables~\guillemotleft~pathologiques\,\footnote{Remarquons qu'il n'existe
aucun objet projectif~\guillemotleft~pathologique~\guillemotright~dans la catégorie $\F$ : tous
les projectifs indécomposables de $\F$ sont de type fini, et  tout
projectif de $\F$ est somme directe de projectifs
indécomposables. Cela provient de la
théorie classique de Krull-Schmidt (cf. \cite{Pop}), parce que les
 projectifs indécomposables de type fini
engendrent $\F$ et ont des anneaux d'endomorphismes locaux.}~\guillemotright~de~$\F$.

 Nous conjecturons notamment le résultat suivant :

\begin{conj}\label{conj-injf} Tout foncteur injectif indécomposable de
  la catégorie $\F$ est à valeurs de dimension finie.
\end{conj}

Les difficultés inhérentes à ce problème, liées à la compréhension
générique des représentations des groupes linéaires, sont
conceptuellement analogues à celles que l'on rencontre pour la
conjecture artinienne extrêmement forte. Cependant, même en admettant la conjecture artinienne extrêmement forte, l'auteur
ignore comment démontrer la conjecture~\ref{conj-injf}, y compris pour la seule enveloppe
injective de~$\bar{P}$, notée~$I_{\bar{P}}$.

\begin{rem}\begin{enumerate}\item On a immédiatement $I_{\bar{P}}(0)=0$
    et $I_{\bar{P}}(\FF)\simeq {\rm hom}_\F
    (\bar{P},I_{\bar{P}})\simeq\FF$. Les résultats de Powell (cf. \cite{GP5})
    fournissent ${\rm hom}_\F(Q_{(2)},I_{\bar{P}})\simeq\FF$ et
    ${\rm hom}_\F(Q_{(2,1)},I_{\bar{P}})\simeq\FF$, d'où l'on déduit,
    compte-tenu de la filtration de $P_{E_2}$ donnée dans \cite{GP2},
    que $I_{\bar{P}}(E_2)$ est de dimension~$7$.
\item Si la conjecture artinienne extrêmement forte est vraie au rang
  $1$ (i.e. pour $\K_1(\F)/\K_0(\F)$), alors il est facile de décrire l'enveloppe injective $I^{\K,1}_{\bar{P}}$
  de
  $\bar{P}$ dans la catégorie $\K_1(\F)$ : c'est l'image par le foncteur $\omega_1$ de
  l'enveloppe injective $I^{\Gr,1}_{(\FF,\FF)}$ (qui est localement finie) du foncteur constant
  $\FF$ de $\F_{\Gr,1}$. Le foncteur $I^{\K,1}_{\bar{P}}$ ainsi obtenu est l'un
  des deux facteurs indécomposables de $\bar{P}\otimes\bar{I}$, dont
  le scindement est donné par application du foncteur $\omega_1$ au
  scindement $\iota_1(\bar{I})\simeq\kappa_1(\bar{I})\oplus
  I^{\Gr,1}_{(\FF,\FF)}$ déduit de
  l'isomorphisme~(\ref{eq-scii}) du §\,\ref{par-prfgr}. 

Mais le foncteur $I^{\K,1}_{\bar{P}}$ est beaucoup moins gros que
l'enveloppe injective $I_{\bar{P}}$ de $\bar{P}$ dans
$\F$, ne serait-ce que parce que cette dernière contient le foncteur
$\overline{\FF[\Gr]}$ (l'article \cite{art2}, §\,$6.3$, montre en
effet que $\overline{\FF[\Gr]}$ est une
extension essentielle de $\bar{P}$).
\end{enumerate}
\end{rem}

La proposition élémentaire suivante va nous aider à discuter la conjecture~\ref{conj-injf}.

\begin{pr}\label{prel-injdf} Le produit tensoriel de deux objets
  injectifs de la catégorie $\F$ dont l'un est à valeurs de dimension
  finie est injectif.
\end{pr}

\begin{proof} Le produit tensoriel entre un injectif de $\F$ et un
  injectif standard $I_V$ est injectif, car l'adjoint à gauche
  $\Delta_V$ à l'endofoncteur $-\otimes I_V$ de $\F$ est exact. 

Si $J$ est un injectif à valeurs de dimension finie, le foncteur
$-\otimes J$ commute aux produits. Comme tout injectif de $\F$ est
facteur direct d'un produit d'injectifs standard, l'injectivité des
$J\otimes I_V$ établie précédemment donne la conclusion.
\end{proof}

\begin{cor}\label{crel-injdf} Supposons que la conjecture artinienne et la
  conjecture~\ref{conj-injf} sont satisfaites. Alors le produit
  tensoriel de deux injectifs de la catégorie $\F$ est injectif.
\end{cor}

\begin{proof} Si la catégorie $\F$ est localement noethérienne, tout
  injectif de $\F$ est isomorphe à une somme directe d'injectifs
  indécomposables (cf. \cite{Gab}). La conclusion résulte donc de la
  proposition~\ref{prel-injdf} et de la distributivité du produit
  tensoriel par rapport à la somme directe.
\end{proof}

Notons $G_0^{inj}(\F)$ le groupe de Grothendieck des objets injectifs
de $\F$ à valeurs de dimension finie de $\F$. Sous les hypothèses du
corollaire~\ref{crel-injdf}, ce groupe hérite d'une structure d'anneau
commutatif induite par le produit tensoriel. Noter que, comme une somme
directe infinie non triviale d'objets injectifs à valeurs de dimension
finie peut être encore à valeurs de dimension finie, $G_0^{inj}(\F)$
est naturellement un anneau topologique non discret. Une
  autre question naturelle (toujours en admettant la conjecture artinienne et la
  conjecture~\ref{conj-injf}) consiste à savoir si le produit
  tensoriel de deux injectifs indécomposables est une somme directe
  {\em finie} d'injectifs indécomposables, ce qui permettrait
  essentiellement de réduire l'étude de $G_0^{inj}(\F)$ au sous-groupe
discret engendré par les classes des injectifs indécomposables (qui en
serait alors un sous-anneau).

Nous verrons dans la section~\ref{sctfco} que la conjecture artinienne
extrêmement forte donnerait une description assez explicite de l'anneau de
Grothendieck $G_0^{tf}(\F)$ des foncteurs de type fini de $\F$. Même en admettant la conjecture~\ref{conj-injf} et la
conjecture artinienne extrêmement forte, la structure de l'anneau
$G_0^{inj}(\F)$ reste obscure.

\part{Préliminaires relatifs aux foncteurs $\omega$ et $\tilde{\nabla}_n$}\label{ptr}

Nous donnons dans cette partie le substrat technique nécessaire à la
généralisation des arguments de la section~\ref{sct-fom} au foncteur
$\omega_n$, où $n$ est un entier supérieur à $1$. Il s'agit d'une part
d'étendre le lemme~\ref{lm-elm}, en remplaçant le foncteur différence
par le foncteur $\tilde{\nabla}_n$, ce que nous accomplissons dans la
section~\ref{par-nablom}. Nous aurons besoin d'autre part d'un
résultat d'annulation cohomologique relatif aux foncteurs
$\tilde{\nabla}_n$-nilpotents (généralisant le cas des foncteurs
finis) ; c'est l'objet de la section~\ref{sct-nnil}. Les deux facettes
de ce programme reposent sur des arguments élémentaires issus de la structure des
injectifs des catégories de foncteurs en grassmanniennes ; le
lemme~\ref{lm-filtr-pol} relatif à la structure de l'injectif $I$ de
$\F$ constitue pratiquement leur seul
ingrédient non formel.

\begin{conv} Dans toute cette partie, on se donne un entier~$n>0$. 
\end{conv}

\section{Les foncteurs $\nabla^\Gr_n$}

Nous introduisons l'outil interne à la catégorie $\F_\Gr$ qui
permettra, dans les sections suivantes, de mener à bien nos
investigations relatives aux foncteurs $\omega$ et $\tilde{\nabla}_n$.

\begin{defi} On définit un endofoncteur
$\nabla^\Gr_n$ de la catégorie $\F_\Gr$ par
$$\nabla^\Gr_n=ker\,\big((- : \iota(I))\to (- :
\iota(p_{n-1}(I)))\big)$$
où la transformation naturelle est induite par l'inclusion
$p_{n-1}(I)\hookrightarrow I$.
\end{defi}

Le lien avec le foncteur $\tilde{\nabla}_n$ est donné par
la proposition suivante.

\begin{pr}\label{lien-nab} Il existe un
  isomorphisme
  $$\omega\circ\nabla^\Gr_n\simeq\tilde{\nabla}_n\circ\omega$$ de
  foncteurs $\F_\Gr\to\F$.
\end{pr}

\begin{proof} Grâce à la proposition~\ref{domeg} et à l'exactitude du
  foncteur $\omega$, il existe un diagramme commutatif aux lignes
  exactes 
$$\xymatrix{ 0\ar[r] & \omega\circ\nabla^\Gr_n\ar@{-->}[d]^\simeq\ar[r] & \omega\circ (- : \iota(I))\ar[d]^\simeq\ar[r]
  & \omega\circ (- : \iota(p_{n-1}(I)))\ar[d]^\simeq
  \\
0\ar[r] & \tilde{\nabla}_n\circ\omega\ar[r] & (- : I) \circ\omega\ar[r] & ( - : p_{n-1}(I))\circ\omega
}$$
dans lequel les flèches verticales sont des isomorphismes.
\end{proof}

Nous donnons maintenant des propriétés formelles du
foncteur $\nabla^\Gr_n$, analogues à celles du foncteur
$\tilde{\nabla}_n$.

La première d'entre elle est une conséquence directe de la suite
exacte (\ref{eq-presi}) de la section~\ref{par-fpol}. On rappelle que
$t_n^* : I\to I_{E_n}$ désigne le morphisme induit par $t_n\in\FF[{\rm
  hom}(E_n,\FF)]$.

\begin{pr}\label{prnab1} Le foncteur $\nabla^\Gr_n$ est l'image de la
  transformation naturelle $(t_n)_*=(- : \iota(t_n^*)) : (- : \iota(I_{E_n}))\to (- : \iota(I))$.
\end{pr}

Comme les foncteurs $(- : \iota(I_{E_n}))$ et $(- : \iota(I))$ sont
exacts (cf. proposition~\ref{divif}), on en déduit :

\begin{cor}\label{crnab-epi} Le foncteur $\nabla^\Gr_n$ est additif ; il préserve les
  monomorphismes et les épimorphismes.
\end{cor}

Nous utiliserons souvent le foncteur $\nabla^\Gr_n$ sur des foncteurs
appartenant à une sous-catégorie $\F_{\Gr,k}$ de $\F_\Gr$. Pour
ramener son étude à  des considérations internes à la catégorie
$\F_{\Gr,k}$, nous emploierons la proposition suivante.

\begin{pr}\label{lm-respn} Soient $k\in\mathbb{N}$, $X$ un objet de
  $\F_{\Gr,k}$ et $\tilde{X}$ l'objet de $\F_\Gr$ image de $X$ par le
  foncteur de prolongement par zéro. Il existe un isomorphisme naturel
  entre la restriction de $\nabla^\Gr_n(\tilde{X})$ à $\F_{\Gr,k}$ et
  l'image de $(t_n)_*=(X : \iota_k(t_n^*)) : (X : \iota_k(I_{E_n}))\to
  (X : \iota_k(I))$.

De plus, le niveau du foncteur $\nabla^\Gr_n(\tilde{X})$ est
inférieur ou égal à $k$.
\end{pr}

\begin{proof} Cela résulte de la proposition~\ref{prnab1} et de l'isomorphisme naturel
  $\mathcal{R}_{\leq k,k}(\mathcal{P}_{k,\leq k}(X):\iota_{\leq
    k}(F))\simeq (X : \iota_k(F)$ (pour $F\in {\rm Ob}\,\F$) de la
  proposition $9.2$ de \cite{art2}, et de l'isomorphisme naturel
  $(\tilde{X} : \iota(F))\simeq\mathcal{P}_{\leq k,\mathbb{N}}(\mathcal{P}_{k,\leq k}(X):\iota_{\leq
    k}(F))$ qui
  s'établit par un argument formel d'adjonction analogue.
\end{proof}

La proposition suivante découle de la commutativité du produit
tensoriel de $\F_\Gr$ et de ce que ses foncteurs de décalage et son
foncteur différence sont isomorphes à des foncteurs de division.

\begin{pr}\label{prnab2} Le foncteur $\nabla^\Gr_n$ commute aux
  foncteurs de décalage et au foncteur différence de $\F_\Gr$, à
  isomorphisme naturel près.
\end{pr}

Nous en venons maintenant à des propriétés élémentaires du foncteur
$\nabla^\Gr_n$ qui diffèrent du comportement de
$\tilde{\nabla}_n$. Elles reposent sur la proposition suivante.

\begin{pr}\label{nab-injp} Il existe dans $\F_\Gr$ un diagramme
  commutatif
\begin{equation}\label{dcf1}\xymatrix{\FF\ar[r]\ar[d] &
  \iota(I)\ar[d]^{\iota(t_n^*)} \\
I^\Gr_{(E_n,E_n)}\ar@{^{(}->}[r] & \iota(I_{E_n})
}\end{equation}
qui induit dans $\F_{\Gr,n}$ un diagramme commutatif  
\begin{equation}\label{dcf2}\xymatrix{\FF\ar[r]\ar@{^{(}->}[d] &
  \iota_n(I)\ar[d]^{\iota_n(t_n^*)} \\
I^{\Gr,n}_{(E_n,E_n)}\ar@{^{(}->}[r] & \iota_n(I_{E_n}).
}\end{equation}
\end{pr}

\begin{proof} Le morphisme $\FF\to I^\Gr_{(E_n,E_n)}$ du
  diagramme~(\ref{dcf1}) est l'unique flèche non nulle ; le morphisme
  $\FF\to\iota(I)$ est la composée de l'unique morphisme non nul
  $\FF\to I^\Gr_{(\FF,\FF)}$ et du monomorphisme scindé
  $I^\Gr_{(\FF,\FF)}\hookrightarrow\iota(I)$ donné par
  l'isomorphisme~(\ref{eq-scii}) de la section~\ref{par-prfgr} ; le
  monomorphisme non spécifié provient de cet isomorphisme également. 

Comme $\FF\simeq\rho(\FF)$, par adjonction entre les foncteurs $\rho$
et $\varepsilon$ (proposition~\ref{divrho}), la commutation du
diagramme~(\ref{dcf1}) se ramène au lemme~\ref{lmsrj-n} ci-dessous,
puisque $\varepsilon\circ\iota\simeq o$ et
$\varepsilon(I^\Gr_{(V,V)})\simeq I^{surj}_V$. 

Dans le diagramme~(\ref{dcf2}) qu'on déduit de~(\ref{dcf1}) par
application du foncteur de restriction, la flèche verticale $\FF\to
I^{\Gr,n}_{(E_n,E_n)}$ est l'unique morphisme non nul, qui est
injectif parce que le foncteur constant $\FF$ de $\F_{\Gr,n}$ est simple. 
\end{proof}

\begin{lm}\label{lmsrj-n}
Le diagramme suivant de $\F_{surj}$ commute
$$\xymatrix{\FF\ar[r]\ar[d] & I^{surj}_\FF\ar@{^{(}->}[r] &
  o(I)\ar[d]^-{o(t^*_n)} \\
I^{surj}_{E_n}\ar@{^{(}->}[rr] & & o(I_{E_n})
}$$
où les deux flèches de source $\FF$ sont les uniques morphismes non
nuls et les deux injections non spécifiées les monomorphismes scindés
donnés par l'isomorphisme canonique $o(I_V)\simeq\bigoplus_{W\in\Gr(V)}I^{surj}_W$.
\end{lm}

\begin{proof}Soit $W$ un sous-espace de $E_n$. Le morphisme
$$\FF\to I^{surj}_\FF\hookrightarrow o(I)\xrightarrow{o(t^*_n)}
o(I_{E_n})\twoheadrightarrow I^{surj}_W$$
est la somme sur les formes linéaires $l\in E_n^*$ telles que
$l(W)\neq 0$ des composées $$\FF\to
I^{surj}_\FF\xrightarrow{l^*}I^{surj}_W,$$  qui sont toutes
égales à l'unique morphisme non nul. Comme le cardinal de l'ensemble $\{l\in
E^*_n\,|\,l(W)\neq 0\}=E_n^*\setminus W^\perp$ est $2^n-2^{n-\dim W}$,
qui est impair si et seulement si $\dim W=n$, on en
conclut que le morphisme $\FF\to I^{surj}_W$ que l'on étudie est non nul
si et seulement si $W=E_n$, d'où l'on déduit le lemme.
\end{proof}

\begin{nota} On désigne par $q^\nabla_n : \nabla^\Gr_n\to id$ la transformation naturelle composée
  de l'inclusion canonique $\nabla^\Gr_n\hookrightarrow ( - : \iota(I))$
  et de la transformation naturelle $( - : \iota(I))\to (- :
  \FF)\simeq id$ induite par le morphisme $\FF\to\iota(I)$ du
  diagramme~(\ref{dcf1}).
\end{nota}


Ainsi, les propositions~\ref{nab-injp} et~\ref{prnab1} procurent un
diagramme commutatif
$$\xymatrix{\nabla^\Gr_n\ar@{^{(}->}[dd]\ar@/^3pc/[ddrr]^-{q^\nabla_n} & & \\
&   (-:\iota(I_{E_n}))\ar@{>>}[lu]\ar[dr]\ar[dl]_-{(t_n)_*} \\
(-:\iota(I))\ar[rr] &  & id.
}$$

\section{Estimation de $\tilde{\nabla}_n\,\omega_n$}\label{par-nablom}

Il n'est pas difficile de voir que, pour $i<n$, le foncteur
$\tilde{\nabla}_i$ a tendance à faire exploser la taille des foncteurs
du type $\omega_n(X)$. La proposition suivante montre que le
comportement de $\tilde{\nabla}_i$ sur $\omega_n(X)$ est tout-à-fait différent pour
$i>n$, lorsque $X$ est un foncteur fini de~$\F_{\Gr,n}$.

\begin{pr}\label{pr-nno} Si $F$ est un foncteur oméga-adapté de
  hauteur strictement inférieure à $n$, alors $F$ est $\tilde{\nabla}_n$-nilpotent.
\end{pr}

\begin{proof} Si $X\in {\rm Ob}\,\F_{\Gr,k}$ est simple, il existe un
  objet simple $F$ de $\F$ et un $GL_k$-module simple $S$ tels que
  $X\simeq\kappa_k(F)\otimes\rho_k(S)$ (cf. proposition~\ref{prcqd}), d'où un épimorphisme
$$P_{E_k}\otimes F\twoheadrightarrow\omega_k\rho_k(S)\otimes F\simeq\omega_k(\iota_k(F)\otimes\rho_k(S))\twoheadrightarrow\omega_k(X).$$
Les propositions~\ref{nab-proj} et~\ref{nab-tens} montrent que
$P_{E_k}\otimes F$, donc $\omega_k(X)$, est
  $\tilde{\nabla}_n$-nilpotent si $k<n$. On conclut grâce à la proposition~\ref{nabnil-ep}.
\end{proof}

Nous nous attachons, dans la suite de cette section, à montrer que
$\tilde{\nabla}_n\,\omega_n(X)$ et $\omega_n(X)$ sont en quelque sorte
du même ordre de grandeur pour $X$ fini. Dans le cas où $X$ est un
foncteur simple pseudo-constant, Powell a calculé exactement, dans \cite{GP2},
$\tilde{\nabla}_n\,\omega_n(X)$ (qui est alors isomorphe à $\omega_n(X)$,
sauf dans un cas) ; nous n'en aurons pas usage.

\begin{nota} Dans cette section, on note $\tilde{X}$ le prolongement
  par zéro à $\F_\Gr$ d'un foncteur $X$ de $\F_{\Gr,n}$. On désigne
  également par $q^\nabla_{n,X}$ le morphisme naturel
  $(q^\nabla_n)_{\tilde{X}} : \nabla^\Gr_n(\tilde{X})\to\tilde{X}$.
\end{nota}

On remarque que la proposition~\ref{lien-nab} entraîne l'existence
d'un isomorphisme naturel
\begin{equation}\label{eq-lnab} \tilde{\nabla}_n\,\omega_n(X)\simeq\omega\nabla^\Gr_n(\tilde{X}).
\end{equation}

En appliquant le foncteur $\omega$ au morphisme naturel
$q^\nabla_{n,X}$, on en déduit un morphisme naturel
\begin{equation}\label{eq-tnnab} \pi_{n,X} : \tilde{\nabla}_n\,\omega_n(X)\to\omega_n(X).
\end{equation}

Pour $n=1$, on retrouve le morphisme $\pi_X$ de la
section~\ref{sct-fom}.

Une description explicite du morphisme $\pi_{n,X}$ sera donnée dans la
section~\ref{par-tsg} (lemme~\ref{pronab}).

\begin{pr} Pour tout foncteur $X$ de $\F_{\Gr,n}$, le morphisme
  $(q^\nabla_n)_{\tilde{X}} : \nabla^\Gr_n(\tilde{X})\to\tilde{X}$ est
  surjectif.
\end{pr}

\begin{proof} Par la proposition~\ref{lm-respn}, cela découle de la surjectivité de la flèche en pointillé du diagramme
$$\xymatrix{X & (X:\iota_n(I))\ar[l] \\
(X:I^{\Gr,n}_{(E_n,E_n)})\ar@{>>}[u] & (X:\iota_n(I_{E_n}))\ar@{-->}[lu]\ar@{>>}[l]\ar[u]_{(t_n)_*}
}$$
induit par le diagramme~(\ref{dcf2}) de la proposition~\ref{nab-injp}. 
\end{proof}

\begin{cor}\label{crf-nabp} La transformation naturelle $\pi_{n} :
  \tilde{\nabla}_n\,\omega_n\to\omega_n$ est surjective.
\end{cor}

\begin{lm}\label{lm1-ngr} Si $X$ est un foncteur pseudo-constant de
  $\F_{\Gr,n}$, la restriction à $\F_{\Gr,n}$ du morphisme
  $(q^\nabla_n)_{\tilde{X}} : \nabla^\Gr_n(\tilde{X})\to\tilde{X}$ est
  un isomorphisme.
\end{lm}

\begin{proof} La proposition $9.9$ de \cite{art2} fournit des
  isomorphismes $(\rho_n(M) : \iota(F))\simeq\rho_n(M :
  F(E_n))\simeq\rho_n(M\otimes F(E_n)^*)$ naturels en le $GL_n$-module
  $M$ et le foncteur $F$ de $\F^{df}$. On en déduit, par la
  proposition~\ref{lm-respn}, des isomorphismes naturels
\begin{eqnarray*}\mathcal{R}_{\mathbb{N},n}\nabla^\Gr_n(\tilde{X}) &
  \simeq & ker\,\big(\rho_n(M\otimes
I(E_n)^*)\xrightarrow{(t_n)_*}\rho_n(M\otimes
p_{n-1}(I)(E_n)^*)\big)\\
& \simeq & \rho_n(M\otimes (I/p_{n-1}(I))(E_n))^*)
\end{eqnarray*}
où $X=\rho_n(M)$ est un foncteur pseudo-constant de $\F_{\Gr,n}$. La
conclusion provient alors de ce que $(I/p_{n-1}(I))(E_n)=\FF$ (lemme~\ref{lm-filtr-pol}).
\end{proof}

\begin{pr}\label{pr2-ngr} Soit $X$ un objet fini de degré $i\geq 0$ de $\F_{\Gr,n}$. Il
  existe dans $\F_\Gr$ une suite exacte $0\to Y\to ker\,q^\nabla_{n,X}\to Z\to 0$
  dans laquelle :
\begin{itemize}\item $Y$ est un foncteur fini de niveau strictement
  inférieur à $n$ et de degré au plus~$i$ ;
\item $Z$ est un foncteur fini de degré strictement inférieur à $i$, de niveau au plus $n$ et de coniveau au moins~$n$. 
\end{itemize}
\end{pr}

\begin{proof} On pose $Y=\mathcal{P}_{\leq
    n-1,\mathbb{N}}\,\mathcal{R}_{\mathbb{N},\leq
    n-1}\nabla^\Gr_n(\tilde{X})$ --- c'est un sous-foncteur de
  $ker\,q^\nabla_{n,X}$ puisque la restriction à $\F_{\Gr,\leq n-1}$
  de $\tilde{X}$ est nulle ---  et $Z=(ker\,q^\nabla_{n,X})/Y$. La
  proposition~\ref{lm-respn} montre que la condition de niveau sur $Z$
  est vérifiée ; la condition de niveau sur $Y$ et celle de coniveau
  sur $Z$ le sont pour des raisons formelles.

La proposition~\ref{prnab2} montre que le foncteur $ker\,q^\nabla_n : \F_{\Gr,n}\to\F_\Gr$
commute au foncteur différence à isomorphisme naturel près. On en
déduit en particulier $\deg Y\leq\deg ker\,q^\nabla_{n,X}\leq\deg
X=i$, par récurrence sur $i$. L'inégalité $\deg Z<i$ s'obtient de
même, par récurrence sur $i$, à partir du lemme~\ref{lm1-ngr}. Comme
les foncteurs $Y$ et $Z$ sont, comme $\nabla^\Gr_n(X)$, à valeurs de
dimension finie, ce qui précède montre qu'ils sont finis, ce qui
achève la démonstration.
\end{proof}

\begin{nota} Dans la proposition suivante,  $\F_{\tilde{\nabla}_n,<i}$
  désigne la plus petite sous-catégorie
  épaisse de $\F$ contenant $\N il_{\tilde{\nabla}_n}$ et les objets
$\omega_n(Y)$, pour $Y$ objet fini de $\F_{\Gr,n}$ de degré
strictement inférieur à $i$.
\end{nota}

Intuitivement, cette catégorie contient tous les objets qui sont~\guillemotleft~moins gros~\guillemotright~que les foncteurs du type $\omega_n(A)$, où $A\in {\rm
  Ob}\,\F_{\Gr,n}$ est de degré~$i$ (la proposition~\ref{pr-nno}
montre que $\F_{\tilde{\nabla}_n,<i}$ contient déjà tous les objets du
type $\omega_k(X)$ avec $k<n$ et $X\in {\rm Ob}\,\F_{\Gr,k}$ fini).

\begin{cor}\label{kerpi-cr} Soient $i\in\mathbb{N}$ et $X$ un objet fini
  de degré $i$ de $\F_{\Gr,n}$.

Le noyau de $\pi_{n,X} : \tilde{\nabla}_n\,\omega_n(X)\to\omega_n(X)$ appartient à $\F_{\tilde{\nabla}_n,<i}$.
\end{cor}

\begin{proof} Soit $0\to Y\to ker\,q^\nabla_{n,X}\to Z\to 0$ la suite
  exacte donnée par la proposition~\ref{pr2-ngr} : par exactitude du
  foncteur $\omega$, on en déduit une suite exacte
  $0\to\omega(Y)\to ker\,\pi_{n,X}\to\omega(Z)\to 0$. Le foncteur
  $\omega(Y)$ est $\tilde{\nabla}_n$-nilpotent, donc objet de
  $\F_{\tilde{\nabla}_n,<i}$, par la proposition~\ref{pr-nno}, et
  $\omega(Z)$ est objet  de
  $\F_{\tilde{\nabla}_n,<i}$ parce que $\deg Z<i$, d'où la proposition.
\end{proof}

\section{Théorème d'annulation cohomologique relatif à la
  $\nabla$-nilpotence}\label{sct-nnil}

La section précédente précisait le comportement du foncteur $\tilde{\nabla}_n$
sur l'image du foncteur $\omega_n$ en terme
de~\guillemotleft~taille~\guillemotright~des objets ; nous allons donner
l'aspect cohomologique de la comparaison entre $\tilde{\nabla}_n$ et $\omega_n$.

\smallskip

Tensorisons le diagramme (\ref{dcf2}) de la proposition~\ref{nab-injp} par un
objet $X$ de $\F_{\Gr,n}$ et appliquons le foncteur $\omega_n$ :  on
obtient un diagramme commutatif naturel
\begin{equation}\label{ome-nd}\xymatrix{\omega_n(X) \ar[r]^-{a_{n,X}}\ar@{^{(}->}[dr]_-{b_{n,X}} & \omega_n(X)\otimes
  I\ar[d]^{\omega_n(X)\otimes t_n^*} \\
& \omega_n(X)\otimes I_{E_n}.
}\end{equation}

Par adjonction entre les endofoncteurs $\Delta_V$ et $-\otimes I_V$ de
$\F$, on en déduit un diagramme commutatif
\begin{equation}\label{eq-omnab}\xymatrix{\Delta_\FF\omega_n(X)\ar[dr]^-{\alpha_{n,X}} & \\
\Delta_{E_n}\omega_n(X)\ar[r]_-{\beta_{n,X}}\ar[u]^{(t_n)_*} & \omega_n(X).
}\end{equation} 

Le théorème fondamental suivant constitue une variation sur la
proposition~\ref{fondcr} adaptée au contexte de la
$\tilde{\nabla}_n$-nilpotence. Nous y reviendrons à la remarque~\ref{lremqnsar}.

\begin{theo}\label{nabpft} Pour tout objet $X$ de $\F_{\Gr,n}$, le
  foncteur $\omega_n(X)$ est $\overline{\N
    il}_{\tilde{\nabla}_n}$-parfait : pour tout foncteur $F$
  de $\overline{\N
    il}_{\tilde{\nabla}_n}$, on a ${\rm Ext}_\F^*(F,\omega_n(X))=0$.
\end{theo}

\begin{proof} Un argument de colimite montre qu'il suffit d'établir la
  nullité des groupes d'extensions ${\rm Ext}^i_\F (F,\omega_n(X))$
  lorsque $F$ est $\tilde{\nabla}_n$-nilpotent, assertion que
  l'on établit par récurrence sur~$i$.

La naturalité de la transformation $(t_n)_* : \Delta_{E_n}\to\Delta_\FF$ se traduit par la
commutation des diagrammes
$$\xymatrix{{\rm hom}_\F
  (F,G)\ar[r]^-{(\Delta_{E_n})_*}\ar[d]_{(\Delta_\FF)_*} & {\rm
    hom}_\F (\Delta_{E_n} F,\Delta_{E_n} G)\ar[d]^{t_n(G)_*} \\
{\rm hom}_\F (\Delta_\FF F,\Delta_\FF G)\ar[r]^-{t_n(F)^*} & {\rm
    hom}_\F (\Delta_{E_n} F,\Delta_\FF G)
}$$
pour tous objets $F$ et $G$ de $\F$, où l'on a noté $t_n(F)$ pour
$((t_n)_*)_F$, pour simplifier les écritures. Comme les foncteurs
$\Delta_{E_n}$ et $\Delta_\FF$ sont exacts, ce diagramme s'étend en un
diagramme commutatif
$$\xymatrix{{\rm Ext}^i_\F
  (F,G)\ar[r]^-{(\Delta_{E_n})_*}\ar[d]_{(\Delta_\FF)_*} & {\rm
    Ext}^i_\F (\Delta_{E_n} F,\Delta_{E_n} G)\ar[d]^{t_n(G)_*} \\
{\rm Ext}^i_\F (\Delta_\FF F,\Delta_\FF G)\ar[r]^-{t_n(F)^*} & {\rm
    Ext}^i_\F (\Delta_{E_n} F,\Delta_\FF G)
}$$
pour tout $i\in\mathbb{N}$.

On forme alors le diagramme commutatif
\begin{equation}\label{ldqt}\xymatrix{{\rm Ext}^i_\F
  (F,\omega_n(X))\ar[r]^-{(\Delta_{E_n})_*}\ar[d]_{(\Delta_\FF)_*} & {\rm
    Ext}^i_\F (\Delta_{E_n} F,\Delta_{E_n}
  \omega_n(X))\ar[d]^{t_n(\omega_n(X))_*}\ar[r]^-\simeq & {\rm
    Ext}^i_\F (F,\Delta_{E_n}\omega_n(X)\otimes
  I_{E_n})\ar[d]^{(t_n(\omega_n(X))\otimes I_{E_n})_*} \\
{\rm Ext}^i_\F (\Delta_\FF F,\Delta_\FF \omega_n(X))\ar[r]^-{t_n(F)^*}\ar[d]_-{(pl^{\nabla,n}_F)^*} & {\rm
    Ext}^i_\F (\Delta_{E_n} F,\Delta_\FF \omega_n(X))\ar[d]^-{(\alpha_{n,X})_*}\ar[r]^-\simeq & {\rm Ext}^i_\F (F,\Delta_\FF\omega_n(X)\otimes I_{E_n})\ar[d]^-{(\alpha_{n,X}\otimes I_{E_n})_*} \\
{\rm Ext}^i_\F (\tilde{\nabla}_n F,\Delta_\FF \omega_n(X))\ar[ru]_-{(pr^{\nabla,n}_F)^*}\ar[d]_-{(\alpha_{n,X})_*} & {\rm
    Ext}^i_\F (\Delta_{E_n} F,\omega_n(X))\ar[r]^-\simeq & {\rm Ext}^i_\F (F,\omega_n(X)\otimes I_{E_n})\\
{\rm Ext}^i_\F (\tilde{\nabla}_n F,\omega_n(X))\ar[ru]_-{(pr^{\nabla,n}_F)^*} & 
}\end{equation}
où $pr^{\nabla,n}_F$ désigne la projection canonique $\Delta_{E_n}
F\twoheadrightarrow\tilde{\nabla}_n F$, et $pl^{\nabla,n}_F :
\tilde{\nabla}_n F\hookrightarrow\Delta_\FF F$ l'inclusion canonique,
et où les isomorphismes de droite dérivent de l'adjonction entre les
endofoncteurs exacts $\Delta_V$ et $-\otimes I_V$ de
$\F$.

Le morphisme $u : {\rm Ext}^i_\F
  (F,\omega_n(X))\to  {\rm
    Ext}^i_\F (F,\omega_n(X)\otimes I_{E_n})$ obtenu en suivant le diagramme est induit
  par le monomorphisme $b_{n,X} :
  \omega_n(X)\hookrightarrow\omega_n(X)\otimes I_{E_n}$. En effet, la
  composée horizontale supérieure est induite par l'unité $\omega_n(X)\to\Delta_{E_n}\omega_n(X)\otimes
  I_{E_n}$ de l'adjonction, tandis que la composée verticale de droite
  est induite par $\beta_{n,X}\otimes I_{E_n}$ (par la commutation du diagramme
  (\ref{eq-omnab})) ; on conclut par adjonction entre les
  morphismes $b_{n,X}$ et $\beta_{n,X}$.

Comme $b_{n,X}$ provient par application du
  foncteur $\omega_n$ à un monomorphisme $X\hookrightarrow
  X\otimes\iota_n(I_{E_n})$, il existe un objet $Y$ de $\F_{\Gr,n}$ et
  une suite exacte courte
$$0\to\omega_n(X)\xrightarrow{b_{n,X}}\omega_n(X)\otimes
I_{E_n}\to\omega_n(Y)\to 0\,.$$
La suite exacte longue de cohomologie associée et l'hypothèse de
récurrence ${\rm Ext}^{i-1}_\F (F,\omega_n(Y))=0$ (on rappelle que $F$
est supposé $\tilde{\nabla}_n$-nilpotent) montrent que le
morphisme  $u : {\rm Ext}^i_\F
  (F,\omega_n(X))\to  {\rm
    Ext}^i_\F (F,\omega_n(X)\otimes I_{E_n})$ est {\em injectif}. On en déduit, en considérant la
  colonne de gauche du diagramme (\ref{ldqt}), un {\em monomorphisme} ${\rm Ext}^i_\F
  (F,\omega_n(X))\hookrightarrow {\rm Ext}^i_\F (\tilde{\nabla}_n
  F,\omega_n(X))$, ce pour tous objets $F$ de $\N
  il_{\tilde{\nabla}_n}$ et $X$ de $\F_{\Gr,n}$. Cela fournit, par
  récurrence sur l'indice de $\tilde{\nabla}_n$-nilpotence de $F$, la
  nullité du groupe d'extensions ${\rm Ext}^i_\F
  (F,\omega_n(X))$, ce qui achève la démonstration.
\end{proof}

\begin{rem}\label{lremqnsar} La proposition \ref{pr-nno} montre que le
  théorème~\ref{nabpft} est une généralisation du
  résultat d'annulation cohomologique de la proposition~\ref{fondcr} (dont
  elle ne fournit toutefois pas l'isomorphisme obtenu pour
  $k=n$). Pour $n=1$, ces deux résultats sont identiques. La
  démonstration du théorème~\ref{nabpft} est une conséquence formelle
  de l'injectivité du morphisme $b_{n,X}$ du diagramme (\ref{ome-nd})
  ; la seule difficulté technique réside dans l'inexactitude du
  foncteur $\tilde{\nabla}_n$, qui oblige à transiter par les
  foncteurs exacts $\Delta_\FF$ et $\Delta_{E_n}$ pour passer aux
  groupes d'extensions. Conceptuellement, cette démonstration procède d'idées très
  voisines de celles employées dans \cite{art2} pour établir la proposition~\ref{fondcr}, reposant
 sur des adjonctions entre foncteurs décalages et produits
  tensoriels convenables.
\end{rem}

\part{Résultats fondamentaux}

Nous nous consacrons désormais aux applications des résultats de la
partie précédente à la structure de la catégorie $\F$. La première
d'entre elle, présentée dans la section~\ref{sctfco}, est le
théorème~\ref{thintg} de l'introduction. Des catégories de foncteurs, elle n'utilise les
estimations de la section~\ref{par-nablom}
que dans le cas des foncteurs pseudo-constants (déjà traité dans
\cite{GP2}), et les propriétés de base des objets finis des catégories
$\F$ et $\F_\Gr$. Un autre ingrédient essentiel provient de la théorie
des représentations : la représentation de Steinberg (modulaire) des
groupes linéaires occupe une place privilégiée, parce que les
facteurs de composition du foncteur de Powell qui lui est associé s'avèrent
bien plus faciles à contrôler à l'aide des foncteurs
$\tilde{\nabla}_n$ que dans le cas général. Outre cette observation
assez directe, qui repose uniquement sur le fait que la partition
régulière associée à la représentation de Steinberg a une longueur
égale à son premier terme, nous emploierons une propriété profonde de
la classe de cette représentation dans l'anneau de Grothendieck de
l'algèbre du groupe linéaire, qui caractérise l'idéal qu'elle
engendre.

Le reste de cette partie consiste en la démonstration du théorème de
simplicité généralisé (théorème~\ref{thint-tsg} de l'introduction) et de ses conséquences sur la conjecture
artinienne extrêmement forte (notamment le
théorème~\ref{thintc}). Grâce aux préliminaires de la
partie~\ref{ptr}, nous pourrons emprunter une voie tout-à-fait
analogue à celle de la partie~\ref{pdx}.

Outre qu'ils reposent sur
l'emploi des foncteurs $\omega$ et $\tilde{\nabla}_n$, les deux volets
principaux de cette dernière partie se rejoignent dans l'observation
suivante : la compréhension fine de la catégorie $\F$ nécessite un
contrôle du comportement générique des représentations des groupes linéaires (ou
symétriques). L'application du théorème de simplicité généralisé à la
structure de $P^{\otimes 2}\otimes F$ (pour $F$ fini) comme la
propriété d'injectivité du morphisme induit par le foncteur $\omega :
\F_\Gr\to\F$ entre groupes de Grothendieck ne peuvent s'établir par
les seules méthodes fonctorielles ; une partie du substrat non formel
sous-jacent à la catégorie $\F$
réside manifestement dans des propriétés intrinsèques des groupes linéaires.

\section{Le morphisme $\omega_* : G^f_0(\F_\Gr)\to  G_0^{tf}(\F)\to\widehat{G}_0^f(\F)$}\label{sctfco}

On rappelle que les notations relatives aux groupes de Grothendieck
ont été données à la fin de l'introduction de cet article, et que les
notations liées aux partitions et facteurs de compositions ont été
introduites au paragraphe~\ref{par-simples}.

\begin{nota} Nous désignerons par $\widehat{G}_0^{f}(\F)$ le groupe
abélien $\mathbb{Z}^\mathfrak{p}$ produit de copies de $\mathbb{Z}$
indexées par les partitions régulières.
\end{nota}

Pour toute partition régulière $\lambda$, la fonction de multiplicité $m_\lambda :
{\rm Ob}\,\F^{tf}\to\mathbb{Z}$ est {\em additive} en ce sens que si $0\to A\to
B\to C\to 0$ est une suite exacte courte de $\F$, alors
$m_\lambda(B)=m_\lambda(A)+m_\lambda(C)$. Elle induit donc un
morphisme de groupes, encore noté $m_\lambda$ par abus, de
$G_0^{tf}(\F)$ vers $\mathbb{Z}$.

\begin{nota} Nous noterons, dans cette section, $j_G :
  G_0^{tf}(\F)\to\widehat{G}_0^{f}(\F)$ le morphisme de groupes dont
  les composantes sont données par les $m_\lambda$
  ($\lambda\in\mathfrak{p}$). Nous l'appellerons {\em morphisme canonique}.
\end{nota}

\begin{rem} La composée $\mathbb{Z}[\mathfrak{p}]\simeq G_0^f(\F)\to
  G_0^{tf}(\F)\xrightarrow{j_G}\widehat{G}_0^{f}(\F)=\mathbb{Z}^\mathfrak{p}$ (où le premier
  morphisme est induit par l'inclusion ${\rm
    Ob}\,\F^{f}\hookrightarrow {\rm Ob}\,\F^{tf}$) coïncide avec
  l'inclusion canonique. Via cette inclusion, on peut voir le groupe
  $\widehat{G}_0^{f}(\F)$ comme un complété convenable de $G_0^f(\F)$.
\end{rem}

Le foncteur $\omega : \F_\Gr\to\F$ est exact, respecte la structure
tensorielle, en prenant à la source le produit tensoriel total $\ptt$
(introduit dans \cite{art2}),
et transforme un objet fini de $\F_\Gr$ en un objet de type fini de
$\F$. Par conséquent, il induit un morphisme d'anneaux $\omega_* :
G_0^f(\F_\Gr)\to G_0^{tf}(\F)$. Nous noterons encore, par abus,  $\omega_* :
G_0^f(\F_\Gr)\to \widehat{G}_0^{f}(\F)$ le morphisme de groupes
composé du morphisme précédent et du
morphisme canonique $j_G : G_0^{tf}(\F)\to
\widehat{G}_0^{f}(\F)$. Notre objectif consiste à démontrer que ce
morphisme $\omega_*$ est injectif.

\smallskip Cette section utilise lourdement certaines propriétés de
  la {\em représentation de Steinberg} de $GL_n$. Avec nos conventions
  d'indexation des représentations simples de $GL_n$, $R_{<n>}$ est la
  représentation de Steinberg, où $<n>$ désigne la partition régulière
  \guillemotleft~triangulaire~\guillemotright~$(n,n-1,\dots,1)$ de $n(n+1)/2$. On pourra se reporter à \cite{Mitch} pour la
démonstration de l'identité entre la définition fonctorielle $R_{<n>}$
et des approches plus classiques de la représentation de Steinberg.

Le lecteur pourra également consulter \cite{Hum} pour
  un survol des propriétés de cette représentation.

\smallskip

La suite de cette section fait usage de la notation $Q_\lambda$
introduite au paragraphe~\ref{par-deffgr}, et du symbole
$\lambda_{+i}$ (où $\lambda$ est une partition) introduit dans la notation~\ref{unoti}.

\begin{lm}\label{lm-crst} Soient $n>0$ et $i\geq 0$ des entiers et
  $\lambda$ une partition régulière telle que $\lambda_1=n$. Si
  $\lambda=<n>$, alors $m_{<n>_{+i}}(Q_\lambda)=1$, sinon $m_{<n>_{+i}}(Q_\lambda)=0$.
\end{lm}

\begin{proof} Notons $f_\lambda(i)=m_{<n>_{+i}}(Q_\lambda)$. Le lemme
  découle des quatre points suivants.
\begin{enumerate}\item La fonction $f_\lambda :
  \mathbb{N}\to\mathbb{N}$ est décroissante pour toute partition
  régulière $\lambda$ telle que $\lambda_1=n$.
\item On a $f_\lambda(0)=1$ si $\lambda=<n>$, $0$ sinon.
\item Pour tout $i\in\mathbb{N}$, $m_{<n>_{+i}}(P_{E_n})>0$.
\item Pour tout $i\in\mathbb{N}$, $m_{<n>_{+i}}(P_{E_n})$ est une
  combinaison linéaire des $f_\lambda(i)$, $\lambda$ parcourant les
  partitions régulières telles que $\lambda_1=n$.
\end{enumerate}

On commence par remarquer que si $\mu$ est une partition régulière
telle que $\mu_1<n$, alors le foncteur de Powell $Q_\mu$, qui est un
quotient de $P^{\otimes \mu_1}$, n'a pas de facteur de composition
$S_\nu$, si $\nu$ est une partition régulière de longueur $n$.

Comme le noyau de l'épimorphisme (corollaire~\ref{crf-nabp}) $\pi_{n,\rho_n(R_\lambda)} :
\tilde{\nabla}_n Q_\lambda\twoheadrightarrow Q_\lambda$ de la section
\ref{par-nablom} admet une filtration finie dont les quotients sont des
foncteurs de ce type (cf. proposition~\ref{pr2-ngr}), on a
$f_\lambda(i)=m_{<n>_{+i}}(\tilde{\nabla}_n Q_\lambda)$. Mais
$m_{<n>_{+i}}(\tilde{\nabla}_n Q_\lambda)\geq m_{<n>_{+i+1}}(Q_\lambda)$
par les propositions \ref{nab-simp} et \ref{evid-nab}.\,$1$, d'où
le premier point.

Le second vient des égalités
$$f_\lambda(0)=\dim {\rm hom}_\F (Q_{<n>},Q_\lambda)=\dim {\rm
  hom}_{GL_n} (R_{<n>},R_\lambda)$$
(cf. \cite{GP2}).

Quant au troisième, c'est une conséquence du théorème \ref{fclambda}.

Le dernier provient de ce que $P_{E_n}\simeq\varpi(P^{surj}_{E_n})$ admet une filtration finie dont
  les quotients sont des foncteurs de Powell $Q_\mu$ avec $\mu_1\leq
  n$ et de la remarque précédente sur les $Q_\mu$ pour $\mu_1<n$.
\end{proof}

\begin{lm}\label{lmcr-st} Pour tout $n\in\mathbb{N}$, l'endomorphisme du groupe de
  Grothendieck $G_0^f(_{\FF[GL_n]}\mathbf{Mod})$ induit par le produit
  tensoriel par la représentation de Steinberg $R_{<n>}$ est injectif.
\end{lm}

\begin{proof} D'après un théorème de Ballard et Lusztig, l'idéal de l'anneau $G_0^f(_{\FF[GL_n]}\mathbf{Mod})$
  engendré par $R_{<n>}$ est égal à l'image canonique de
  $K_0(_{\FF[GL_n]}\mathbf{Mod})$ dans
  $G_0^f(_{\FF[GL_n]}\mathbf{Mod})$. Pour une démonstration de ce
  résultat, voir \cite{CR2}, chapitre~$8$, théorème~$72.10$ (la
  proposition~$9.3$ de cet ouvrage montre
  que les groupes linéaires sur $\FF$ vérifient les hypothèses dudit
  théorème). On en déduit, par la théorie
  générale des représentations modulaires (cf. théorème
  $21.22$ de \cite{CR}), que le conoyau de l'endomorphisme de $G_0^f(_{\FF[GL_n]}\mathbf{Mod})$ induit par le produit
  tensoriel par $R_{<n>}$ est {\em fini}. Comme
  $G_0^f(_{\FF[GL_n]}\mathbf{Mod})$ est un $\mathbb{Z}$-module libre
  de rang fini, cet endomorphisme est injectif.
\end{proof}

\begin{nota}\label{tronqp} \'Etant donnés des entiers $n$ et $i$, nous noterons
  $\mathfrak{p}_{n,i}=\{\lambda\in\mathfrak{p}\,|\,\lambda_n\geq i\}$
  (où l'on convient que $\lambda_n=+\infty$ si $n\leq 0$),
  et $\widehat{G}_0^{f}(\F)_{n,i}=\mathbb{Z}^{\mathfrak{p}_{n,i}}$. C'est donc un
  quotient du groupe $\widehat{G}_0^{f}(\F)$.

Si $d$ est un autre entier, nous poserons $\mathfrak{p}_{n,i,\leq
  d}=\{\lambda\in\mathfrak{p}_{n,i}\,|\,|\lambda|\leq d\}$ et
$\widehat{G}_0^{f}(\F)_{n,i, \leq d}=\mathbb{Z}^{\mathfrak{p}_{n,i,
    \leq d}}$.

 Nous noterons enfin $\mathfrak{p}_{<
  i}=\{\lambda\in\mathfrak{p}\,|\,\lambda_1<i\}$, $\mathfrak{p}^{\leq
  n}_{n,i}=\{\lambda\in\mathfrak{p}_{n,i}\,|\,l(\lambda)=n\}$ et $\mathfrak{p}^{\leq
  n}_{n,i,\leq d}=\mathfrak{p}^{\leq  n}_{n,i}\cap\mathfrak{p}_{n,i,\leq  d}$.
\end{nota}

\begin{lm}\label{lmoub} Pour tout $n\in\mathbb{N}$, le morphisme
$$\alpha : \mathbb{Z}[\mathfrak{p}_{<n+1}]\to
G_0^f(_{\FF[GL_n]}\mathbf{Mod})\qquad [\lambda]\mapsto [I_\lambda(E_n)]$$
a un conoyau fini.
\end{lm}

\begin{proof} Notons
  $M=(M_{\lambda,\mu})_{\lambda,\mu\in\mathfrak{p}_{<n+1}}$ la matrice définie par le
  fait que $I_\lambda$ possède une filtration finie dont les
  sous-quotients sont les duaux $J_\mu$ de $Q_\mu$, chacun apparaissant avec la multiplicité
  $M_{\lambda,\mu}$. L'analyse effectuée dans le §\,$4$ de l'article
  \cite{GP2} (où la matrice $M$ est notée $(\alpha_{\lambda,\mu})$)
  montre que $M$ est une matrice inversible sur $\mathbb{Q}$, car il
  en est de même pour les matrices de Cartan des groupes finis $GL_i$.

Les colonnes de la matrice
$M^{-1}$ ont pour images par $\alpha\otimes\mathbb{Q}$ les images des
$J_\lambda(E_n)$ (où $\lambda\in\mathfrak{p}_{<n+1}$) dans
$G_0^f(_{\FF[GL_n]}\mathbf{Mod})\otimes\mathbb{Q}$. Comme
$J_\lambda(E_n)\simeq R_\lambda$ si $\lambda_1=n$, cela montre que le
morphisme $\alpha\otimes\mathbb{Q}$ est surjectif, d'où la conclusion.
\end{proof}

\begin{nota} Dans la suite de cette section, si $\alpha$ est une
  partition régulière telle que $\alpha_1=n$, nous noterons
  $P^{GL_n}_\alpha$ la couverture projective de $R_\alpha$ dans
  $_{\FF[GL_n]}{\bf Mod}$.

Nous désignerons par $\Delta_\mu$, pour
$\mu\in\mathfrak{p}$, l'endofoncteur exact $(\cdot : I_\mu)$ de
$\F$. C'est un facteur direct de $\Delta_{E_n}$, où $n=\mu_1$.
\end{nota}

\begin{pr}\label{prdst} Le morphisme
$$\beta : G_0^f(_{\FF[GL_n]}\mathbf{Mod})\to\mathbb{Z}^{\mathfrak{p}_{<n+1}}\qquad
[M]\mapsto
(m_{<n>_{+i}}(\Delta_\mu\omega_n\rho_n(M)))_{\mu}$$
est injectif.
\end{pr}

\begin{proof} Par la proposition~\ref{domeg} et la proposition~$9.2$ de
  \cite{art2} (cf. démonstration de la proposition~\ref{lm-respn}),
il existe un épimorphisme  $$\Delta_\mu
Q_\lambda\twoheadrightarrow\omega_n(\rho_n(R_\lambda):\iota_n(I_\mu))$$
dont le noyau est l'image par $\omega_{\leq n-1}$ d'un foncteur fini
pseudo-constant de $\F_{\Gr,\leq n-1}$. On en déduit
(cf. démonstration du lemme \ref{lm-crst}) l'égalité
$$m_{<n>_{+i}}(\Delta_\mu
Q_\lambda)=m_{<n>_{+i}}(\omega_n(\rho_n(R_\lambda):\iota_n(I_\mu))).$$
On a $(\rho_n(R_\lambda):\iota_n(I_\mu))\simeq\rho_n(R_\lambda:I_\mu(E_n))$ par la proposition~\ref{divrho}. Le lemme
\ref{lm-crst} donne alors $$m_{<n>_{+i}}(\Delta_\mu
Q_\lambda)=m_{R_{<n>}}(R_\lambda\otimes
I_\mu(E_n)^*).$$ On a donc, puisque $\FF$ est un corps de
décomposition de $\FF[GL_n]$,
$$m_{<n>_{+i}}(\Delta_\mu
Q_\lambda)=\dim {\rm hom}_{GL_n}(R_{<n>},R_\lambda\otimes
I_\mu(E_n)^*).$$

Le lemme~\ref{lmoub} montre alors que le noyau du morphisme $\beta$
est inclus dans le noyau du morphisme
$$G_0^f(_{\FF[GL_n]}\mathbf{Mod})\to\mathbb{Z}\qquad\dim {\rm hom}_{GL_n}(R_{<n>},R_\lambda\otimes
M^*)$$
pour tout  $GL_n$-module $M$.

Or on a $$\dim {\rm hom}_{GL_n}(R_{<n>},R_\lambda\otimes
(P_\mu^{GL_n})^*)=\dim {\rm hom}_{GL_n}(R_{<n>}\otimes
P_\mu^{GL_n},R_\lambda)$$
$$=\dim {\rm
  hom}_{GL_n}(P_\mu^{GL_n},R_\lambda\otimes
R_{<n>})=m_{R_\mu}(R_\lambda\otimes R_{<n>})\,,$$
où l'on a   utilisé
l'auto-dualité et la projectivité de la représentation de Steinberg
(cf. \cite{Jan}, chapitre $10$, (10.1) et~(10.2)).

L'injectivité du morphisme $\beta$ découle maintenant du lemme~\ref{lmcr-st}
\end{proof}

\begin{lm}\label{prti-g} Soient $n$ et $i$ deux entiers positifs. Il
  existe un entier $d$ tel que le
  morphisme de groupes
$$G_0^f(_{\FF[GL_n]}\mathbf{Mod})\xrightarrow{(\omega_n\rho_n)_*}
G_0^{tf}(\F)\xrightarrow{j_G}\widehat{G}_0^f(\F)\twoheadrightarrow\mathbb{Z}^{\mathfrak{p}_{n,i,
    \leq d}^{\leq n}}$$
est injectif.
\end{lm}

\begin{proof} Soit $\mu\in\mathfrak{p}_{<n+1}$. Pour tout $GL_n$-module fini $M$, on a 
$$m_{<n>_{+i}}(\Delta_\mu\omega_n\rho_n(M))=\sum_{\lambda\in\mathfrak{p}}m_{<n>_{+i}}(\Delta_\mu
S_\lambda)\, m_\lambda(\omega_n\rho_n(M))$$
(somme dont seul un nombre fini des termes sont non nuls). Si
$\lambda\in\mathfrak{p}$ est telle que
$\lambda\vdash\omega_n\rho_n(M)$ et $<n>_{+i}\,\vdash\Delta_\mu
S_\lambda$, on a $l(\lambda)\leq n$ parce que $\omega_n\rho_n(M)$
s'obtient par extensions de quotients de $P_{E_n}$ ; en appliquant
$(\tilde{\nabla}_n)^i$ à la relation $<n>_{+i}\,\vdash\widetilde{\Delta}^n
S_\lambda$, on obtient, compte-tenu des propositions~\ref{evid-nab} et~\ref{nab-simp}, $(\tilde{\nabla}_n)^i S_\lambda\neq
0$, puis $\lambda_n\geq i$. Par ailleurs, comme $l(\lambda)=n$, si $\nu\vdash\widetilde{\Delta}^n
S_\lambda$ et $l(\nu)=n$, alors $|\nu|\geq|\lambda|-n^2$ grâce à la
proposition $4.3.1.2$ de \cite{GP3} et au théorème~\ref{fclambda}, donc
$n(n+1)/2+ni=|<n>_{+i}|\,\geq|\lambda|-n^2$.

Par conséquent, si $d\geq n(n+1)/2+ni+n^2$, on peut compléter le diagramme
$$\xymatrix{G_0^f(_{\FF[GL_n]}\mathbf{Mod})\ar[r]^-\gamma\ar[d]_\beta &
 \mathbb{Z}^{\mathfrak{p}_{n,i,
    \leq d}^{\leq n}}\ar@{-->}[ld] \\
\mathbb{Z}^{\mathfrak{p}_{<n+1}} & 
}$$
où $\beta$ est le morphisme de la
proposition~\ref{prdst}. L'injectivité de $\beta$ donnée par cette
proposition implique donc celle
de $\gamma$.
\end{proof}

\begin{lm}\label{prelpp} Pour tous entiers $n, i\geq 0$, le morphisme de groupes
  $$\mathbb{Z}[\mathfrak{p}_{n,i}^{\leq n}]\otimes\mathbb{Z}[\mathfrak{p}_{<
    i}]\hookrightarrow G_0^{f}(\F)\otimes G_0^{f}(\F)\to G_0^{f}(\F)\twoheadrightarrow\mathbb{Z}[\mathfrak{p}_{n,i}]$$
composé de l'inclusion déduite de l'isomorphisme
$G_0^{f}(\F)\simeq\mathbb{Z}[\mathfrak{p}]$, du produit de l'anneau
$G_0^{f}(\F)$ et de la projection
$G_0^{f}(\F)\simeq\mathbb{Z}[\mathfrak{p}]\twoheadrightarrow\mathbb{Z}[\mathfrak{p}_{n,i}]$
déduite de l'inclusion $\mathfrak{p}_{n,i}\hookrightarrow\mathfrak{p}$, est injectif.
\end{lm}

\begin{proof} Si $\lambda\in\mathfrak{p}_{n,i}^{\leq n}$ et
  $\mu\in\mathfrak{p}_{< i}$, alors la partition $(\lambda,\mu)$ est
  régulière et appartient à $\mathfrak{p}_{n,i}$. De plus, la fonction
  $f : \mathfrak{p}_{n,i}^{\leq n}\times\mathfrak{p}_{< i}\to\mathfrak{p}_{n,i}$ associant aux
  partitions $\lambda$ et $\mu$ la partition $(\lambda,\mu)$ est
  injective, ce qui permet d'identifier cette partition au couple
  $(\lambda,\mu)$. 

Soit $a : \mathfrak{p}_{n,i}^{\leq n}\times\mathfrak{p}_{< i}\to\mathbb{N}$ une bijection telle que
  $a(\alpha)<a(\beta)$ si $|\alpha | <|\beta |$ ou $|\alpha |=|\beta |$ et $\alpha>\beta$, et définissons une matrice
  $(u_{i,j})_{(i,j)\in\mathbb{N}^2}$ par l'égalité
  $u_{a(\lambda,\mu),a(\nu)}=m_\nu(S_\lambda\otimes S_\mu)$ :
c'est la matrice dans la base canonique (indexée
par $a$) du morphisme
$$\mathbb{Z}[\mathfrak{p}_{n,i}^{\leq n}\times\mathfrak{p}_{<
    i}]\simeq\mathbb{Z}[\mathfrak{p}_{n,i}^{\leq n}]\otimes\mathbb{Z}[\mathfrak{p}_{<
    i}]\to\mathbb{Z}[\mathfrak{p}_{n,i}]\twoheadrightarrow\mathbb{Z}[\mathfrak{p}_{n,i}^{\leq n}\times\mathfrak{p}_{<
    i}]$$
où la flèche centrale est celle de l'énoncé et le dernier épimorphisme
est induit par $f$.

Le théorème \ref{fclambda} montre que la matrice
$(u_{i,j})_{(i,j)\in\mathbb{N}^2}$ est triangulaire, et la proposition
\ref{prck} que ses coefficients diagonaux valent $1$, ce qui achève la démonstration.
\end{proof}

\begin{lm}\label{lm-grn} Soient $n$ et $i$ des entiers naturels. Le
  morphisme
$$G_0^f(\F_{\Gr,n})\hookrightarrow
G_0^f(\F_\Gr)\xrightarrow{\omega_*}\widehat{G}^f_0(\F)\twoheadrightarrow\widehat{G}^f_0(\F)_{n,i}$$
est injectif.
\end{lm}

\begin{proof} Soient $j\in\mathbb{N}$ et $G(j)_n$ le sous-groupe de
  $G_0^f(\F_{\Gr,n})$ engendré par les classes d'objets simples de
  degré au plus $j$ ; nous noterons simplement $G(j)$ pour $G(j)_0$. Comme $G_0^f(\F_{\Gr,n})$ est la réunion
  croissante des sous-groupes $G(j)_n$, il suffit de montrer que la
  restriction à $G(j)_n$ du morphisme de l'énoncé est injective pour
  tout~$j$.

La proposition~\ref{prcqd} montre que le foncteur
$\xi_n : \F\otimes\,_{\FF[GL_n]}\mathbf{Mod}\to\F_{\Gr,n}$ induit un isomorphisme
$$G_0^f(_{\FF[GL_n]}\mathbf{Mod})\otimes G(j)\xrightarrow{\simeq}
G(j)_n.$$
Par la proposition \ref{prfig}.\,$4$, on en déduit un diagramme commutatif
$$\xymatrix{G_0^f(_{\FF[GL_n]}\mathbf{Mod})\otimes
G(j)\ar@{^{(}->}[r]\ar[d]_-{(\omega_n\rho_n)_*\otimes G(j)} & G_0^f(\F_{\Gr,n})\ar[d]^-{(\omega_n)_*} \\
G^{tf}_0(\F)\otimes
G(j)\ar[r]\ar[d]_-{j_G\otimes G(j)} & G^{tf}_0(\F)\ar[d]^-{j_G} \\
\widehat{G}^f_0(\F)\otimes G(j)\ar[d] & \widehat{G}^f_0(\F)\ar@{>>}[d] \\
\mathbb{Z}[\mathfrak{p}_{n,i,\leq d}^{\leq n}]\otimes\mathbb{Z}[\mathfrak{p}_{<
    i}]\ar@{^{(}->}[r] & \widehat{G}^f_0(\F)_{n,i}
}$$
pour $i>j$ et un entier $d$ convenable, dans lequel :
\begin{itemize}\item les flèches horizontales sont induites par le
  produit tensoriel --- celle du bas est injective par le lemme
  \ref{prelpp} ;
\item les flèches verticales inférieures sont induites par les
  projections canoniques et l'inclusion $G(j)\subset\mathbb{Z}[\mathfrak{p}_{<i}]$ déduite de l'inégalité $i>j$ (on note que $\mathbb{Z}[\mathfrak{p}_{n,i,\leq d}^{\leq n}]\otimes\mathbb{Z}[\mathfrak{p}_{<i}]\simeq\mathbb{Z}^{\mathfrak{p}_{n,i,\leq d}^{\leq n}}\otimes\mathbb{Z}^{\mathfrak{p}_{<i}}$ par finitude des ensembles  $\mathfrak{p}_{n,i,\leq d}^{\leq n}$ et $\mathfrak{p}_{<i}$) ;
\item l'entier $d$ est choisi, conformément au lemme \ref{prti-g},
  pour que la composée verticale de gauche soit injective.
\end{itemize}

La composée $G_0^f(_{\FF[GL_n]}\mathbf{Mod})\otimes
G(j)\to\widehat{G}^f_0(\F)_{n,i}$ obtenue en suivant le diagramme est
donc injective (suivre la moitié gauche) ; elle s'identifie à la
restriction à $G(j)_n$ du morphisme de l'énoncé (suivre la moitié droite). Celle-ci est donc
injective si $i>j$, donc pour tout $i$, ce qui achève la démonstration. 
\end{proof}

Nous pouvons désormais établir le résultat principal de cette section.

\begin{theo}\label{th-fc} Le morphisme de groupes $\omega_* :
  G_0^f(\F_\Gr)\to\widehat{G}_0^f(\F)$ est injectif.
\end{theo}

\begin{proof} Il suffit d'établir que pour tout $n\in\mathbb{N}$, la
  composée
$$G_0^f(\F_{\Gr,\leq n})\hookrightarrow
G_0^f(\F_\Gr)\xrightarrow{\omega_*}\widehat{G}_0^f(\F)$$
est injective, ou encore que
  pour tout $n$ et tout $j\in\mathbb{N}$, la restriction $f_{n,j}$ de ce
  morphisme au sous-groupe $G(j)_{\leq n}$ de  $G_0^f(\F_{\Gr,\leq
    n})$ engendré par les classes de foncteurs simples de degré au
  plus $j$ est injective. 

Pour $i>j$, la composée $G(j)_{\leq n-1}\hookrightarrow
G_0^f(\F_\Gr)\to\widehat{G}_0^f(\F)\twoheadrightarrow\widehat{G}^f_0(\F)_{n,i}$
est nulle : si $\lambda\vdash P^{\otimes k}\otimes F$ avec $k<n$ et
$\deg F<i$, alors $\lambda\notin\mathfrak{p}_{n,i}$). Le lemme
précédent permet d'en déduire que le noyau de $f_{n,j}$ est inclus
dans  $G(j)_{\leq n-1}$, d'où $ker\,f_{n,j}=ker\,f_{n-1,j}$. On
conclut par récurrence sur $n$.
\end{proof}

\begin{cor}\label{crfcom2}\begin{enumerate}\item Le morphisme d'anneaux $\omega_* :
  G_0^f(\F_\Gr)\to G_0^{tf}(\F)$ est injectif.
\item Si la conjecture artinienne extrêmement forte est vérifiée,
  c'est un isomorphisme, et le morphisme de groupes $j_G  :
  G_0^{tf}(\F)\to\widehat{G}_0^f(\F)$ est injectif.
\end{enumerate}
\end{cor}

\begin{proof} Comme la conjecture artinienne extrêmement forte
  implique que le morphisme $\omega_* :
  G_0^f(\F_\Gr)\to G_0^{tf}(\F)$ est surjectif, la conclusion découle
  du théorème précédent.
\end{proof}

\section{Théorème de simplicité généralisé}\label{par-tsg}

Le but de cette section est d'établir le théorème~\ref{thsg}, qui
constitue un succédané de la conjecture artinienne extrêmement forte,
où la filtration de Krull est remplacée par la filtration par
$\nabla$-nilpotence de~$\F$. Nous suivons la même marche qu'à la
section~\ref{sct-fom}, où les
assertions~\ref{lp3} et~\ref{lp4} du lemme~\ref{lm-elm} sont
remplacées par les corollaires~\ref{crf-nabp} et~\ref{kerpi-cr}
respectivement, et la proposition~\ref{fondcr}.\,$1$ par le théorème~\ref{nabpft}.

\begin{conv} Dans toute cette section, on se donne un entier
  strictement positif~$n$.
\end{conv}

 On rappelle que le morphisme
  $\pi_{n,X} : \tilde{\nabla}_n\omega_n(X)\to\omega_n(X)$ a été défini dans la section~\ref{par-nablom}.

\begin{nota} Si $k$ est un entier naturel et $X$ un
  objet de $\mathcal{F}_{\mathcal{G}r,n}$, nous désignerons par $\pi_{n,X}^k
  : (\tilde{\nabla}_n)^k\omega_n(X)\to\omega_n(X)$ le morphisme
  $\pi_{n,X}\circ\tilde{\nabla}_n\,\pi_{n,X}\circ\dots\circ
  (\tilde{\nabla}_n)^{k-1}\pi_{n,X}$.
\end{nota}

Afin de mener des estimations explicites à l'aide de ces morphismes,
nous identifions le morphisme $\pi_{n,X}$ sur les objets grâce au
lemme suivant, qui emploie le morphisme $\beta_{n,X}$ du
diagramme~(\ref{eq-omnab}) de la section~\ref{par-nablom}.

\begin{lm}\label{pronab} Soit $X\in {\rm Ob}\,\F_{\Gr,n}$.
\begin{enumerate}\item Le morphisme $\beta_{n,X} : \Delta_{E_n}\omega_n(X)\to\omega_n(X)$ est la composée
  de la projection canonique
  $\Delta_{E_n}\omega_n(X)\twoheadrightarrow\tilde{\nabla}_n\,\omega_n(X)$ et du morphisme $\pi_{n,X} : \tilde{\nabla}_n\,\omega_n(X)\to\omega_n(X)$.
\item \'Etant donnés $V\in {\rm Ob}\,\E^f$ et $W\in\Gr_n(V)$, le
  morphisme
$$X(V\oplus E_n,W)\hookrightarrow\omega_n(X)(V\oplus
E_n)=\Delta_{E_n}\omega_n(X)(V)\xrightarrow{(\beta_{n,X})_V}\omega_n(X)(V)$$
est nul si l'une des deux applications linéaires $W\hookrightarrow V\oplus E_n\twoheadrightarrow
E_n$ et $W\hookrightarrow V\oplus E_n\twoheadrightarrow V$ est de rang
$<n$, et est induite par la projection $V\oplus
E_n\twoheadrightarrow V$ sinon.
\end{enumerate}
\end{lm}

\begin{proof} La vérification de la première assertion est
  immédiate. La seconde résulte des observations suivantes, où l'on
  note $\tilde{X}$ le prolongement par zéro à $\F_\Gr$ de~$X$ :
\begin{enumerate}\item le morphisme $\beta_{n,X}$ s'obtient en prenant
  la restriction à $\F_{\Gr,n}$ du morphisme $(\tilde{X} :
  \iota_n(I_{E_n}))\to (\tilde{X} : I^\Gr_{(E_n,E_n)})\to\tilde{X}$ induit
  par le diagramme~(\ref{dcf1}) de la proposition~\ref{nab-injp} puis en appliquant le
  foncteur $\omega_n$ ;
\item par la proposition~\ref{divif}, évalué sur un objet $(V,B)$ de $\E^f_\Gr$, le morphisme  $(\tilde{X} :
  \iota_n(I_{E_n}))\to (\tilde{X} : I^\Gr_{(E_n,E_n)})$ est donné  par la projection
$$\underset{im\,(W\hookrightarrow V\oplus
  E_n\twoheadrightarrow V)=B}{\bigoplus_{W\in\Gr(V\oplus E_n)}} \tilde{X}(V\oplus
E_n,W)\twoheadrightarrow\underset{im\,(W\hookrightarrow V\oplus
  E_n\twoheadrightarrow E_n)=E_n}{\underset{im\,(W\hookrightarrow V\oplus
  E_n\twoheadrightarrow V)=B}{\bigoplus_{W\in\Gr(V\oplus E_n)}}}\tilde{X}(V\oplus
E_n,W) ;$$
\item le morphisme $(\tilde{X} : I^\Gr_{(E_n,E_n)})\to\tilde{X}$,
  évalué sur un objet $(V,B)$, est donné par l'application linéaire
$$\underset{im\,(W\hookrightarrow V\oplus
  E_n\twoheadrightarrow E_n)=E_n}{\underset{im\,(W\hookrightarrow V\oplus
  E_n\twoheadrightarrow V)=B}{\bigoplus_{W\in\Gr(V\oplus E_n)}}}\tilde{X}(V\oplus
E_n,W)\to\tilde{X}(V,B)$$
dont chaque composante est induite par la projection $V\oplus
E_n\twoheadrightarrow V$. Cela provient du lemme~\ref{lm-nabom} et du
lemme de Yoneda.
\end{enumerate}
\end{proof}

\begin{lm}\label{elnab} Soient $V$ un espace vectoriel de dimension finie, $X$ un
  objet de $\mathcal{F}_{\mathcal{G}r,n}$, $H$ un
  sous-espace vectoriel de dimension $n$ de $V^*$ et $x$ un élément
  de $\omega_n(X)(V)$. Notons $a_x : P_V\to \omega_n(X)$ le morphisme
  de $\mathcal{F}$ représenté par $x$, et $s_H=\sum_{h\in H}[h]\in k_n P(V^*)\subset\FF[V^*]$. 

Le morphisme
$$P_V\xrightarrow{t\mapsto s_H\otimes t}k_n P(V^*)\otimes P_V\simeq
\tilde{\nabla}_n\,P_V\xrightarrow{\tilde{\nabla}_n(a_x)}\tilde{\nabla}_n\,\omega_n(X)\xrightarrow{\pi_{n,X}}\omega_n(X)$$
représente l'élément $\Pi^\omega_{X,V,H}(x)$ de $\omega_n(X)(V)$, où
$\Pi^\omega_{X,V,H}$ désigne la projection  de $\omega_n(X)(V)$ sur
les facteurs directs $X(V,W)$, où $W\in\mathcal{G}r_n(V)$ est tel que
$H\cap W^\perp=\{0\}$.
\end{lm}

\begin{proof} Le diagramme
$$\xymatrix{ & \Delta_{E_n} P_V\ar[rr]^-{\Delta_{E_n}a_x}\ar@{>>}[d] & & 
  \Delta_{E_n}\omega_n(X)\ar[dr]^-{\beta_{n,X}}\ar@{>>}[d] & \\
P_V\ar[r]_-{s_H\otimes\cdot}\ar[ru]^-{[(l_1,\dots,l_n)]\otimes\cdot} &
\tilde{\nabla}_n\,P_V\ar[rr]_-{\tilde{\nabla}_n(a_x)} & & 
\tilde{\nabla}_n\,\omega_n(X)\ar[r]_-{\pi_{n,X}} & \omega_n(X)
}$$
commute, où :
\begin{itemize}\item $\{l_1,\dots,l_n\}$ désigne une base de $H$ (on a
  noté $(l_1,\dots,l_n)$ l'élément de ${\rm hom}_\E (V,E_n)$ dont les
  composantes sont les $l_i$) ;
\item l'on a
identifié $\Delta_{E_n} P_V$ à $\FF[{\rm hom}_\E (V,E_n)]\otimes P_V$
(et $\tilde{\nabla}_n\,P_V$ à $k_n P(V^*)\otimes P_V$) ;
\end{itemize}

L'élément $[(l_1,\dots,l_n)]\otimes [id_V]$ de $\Delta_{E_n} P_V(V)$
s'envoie par $\Delta_{E_n}a_x$ sur l'élément de
$\Delta_{E_n}\omega_n(X)(V)=\omega_n(X)(V\oplus E_n)$ image par
$f=(id_V,l_1,\dots,l_n) : V\to V\oplus E_n$ de
$x$.

Soit $W\in\Gr_n(V)$. Le lemme~\ref{pronab} montre que la composée
$$X(V,W)\hookrightarrow\omega_n(X)(V)\xrightarrow{f_*}\omega_n(X)(V\oplus
E_n)\xrightarrow{\beta_{n,X}}\omega_n(X)(V)$$
est nulle si la projection de $f(W)$ sur $E_n$ en est un sous-espace
strict, condition équivalente à la non-inversibilité de la restriction
à $W$ de $(l_1,\dots,l_n)$, ou encore à $H\cap W^\perp\neq\{0\}$, et
que sinon elle est égale à l'inclusion
$X(V,W)\hookrightarrow\omega_n(X)(V)$, puisque la composée de $f :
V\to V\oplus E_n$ et de la projection $V\oplus E_n\twoheadrightarrow V$ est
l'identité. Cela établit le lemme.
\end{proof}

\begin{nota} Dans cette section, si $V$ est un espace vectoriel de
  dimension finie, $W$ un sous-espace de dimension $n$ de $V$ et $H$
  un sous-espace de dimension $n$ de $V^*$, nous noterons $<H,W>$
  l'élément de $\FF$ égal à $1$ si
  $H\cap W^\perp=\{0\}$, $0$ sinon. Autrement dit, si
  $(w_1,\dots,w_n)$ (resp. $(l_1,\dots,l_n)$) est une base de $W$
  (resp. $H$), on a $<H,W>=\det (l_i(w_j))$.
\end{nota}

La propriété de stabilisation suivante, qui fournit la partie~\guillemotleft~concrète~\guillemotright~de la démonstration du théorème~\ref{thsg}, généralise la proposition~\ref{cme1}.

\begin{pr}\label{cme2} Soient $X$ un objet de $\F_{\Gr,n}$ et $F$ un
  sous-objet de $\omega_n(X)$. Pour tout entier $k\geq 0$, on note
  $C_k=\pi_{n,X}^k((\tilde{\nabla}_n)^{k}F)$.
\begin{enumerate}\item La suite $(C_k)_{k>0}$ de sous-objets de $\omega_n(X)$
  est croissante ; nous noterons $C_\infty$ sa réunion. Si $X$ est un objet
  noethérien de $\F_{\Gr,n}$, cette suite stationne. 
\item Pour tout $k\geq 0$, le foncteur $C_k$ est engendré par les éléments
  du type $$(<H_1,W>\dots <H_k,W> x_W)_{W\in\Gr_n(V)}\in\omega_n(X)(V),$$
où $V$ parcourt les espaces vectoriels de
  dimension finie, $x=(x_W)_{W\in\Gr_n(V)}$ les éléments de
  $F(V)$ et $(H_1,\dots,H_k)$ les $k$-uplets d'éléments de $\Gr_n(V^*)$.
\item Le foncteur $C_\infty$ est le plus petit sous-$\bar{G}(n)$-comodule de
  $\omega_n(X)$ contenant $F$. Si $F$ est lui-même un sous-$\bar{G}(n)$-comodule de
  $\omega_n(X)$, on a $C_k=F$ pour tout $k\geq 0$.
\item Si $X$ est localement fini, alors ${\rm hom}\,(C_\infty/F,\omega_n(T))=0$ pour tout objet $T$ de~$\F_{\Gr,n}$.
\end{enumerate}
\end{pr}

\begin{proof} Le lemme \ref{elnab} et la préservation
des épimorphismes par le foncteur $\tilde{\nabla}_n$ montrent le second point
pour $k=1$. Le cas général s'en déduit aussitôt par récurrence.

On en déduit que la suite $(C_k)_{k>0}$ est croissante. Si $W$ est un
sous-espace de dimension $n$ d'un espace vectoriel de dimension finie $V$, soient
$H_1,\dots,H_k$ les éléments de $\Gr_n(V^*)$ tels que $<H_i,W>=1$. Alors la
fonction $<H_1,\cdot>\dots <H_k,\cdot>\: : \,\Gr_n(V)\to\FF$ est égale à
l'indicatrice de~$\{W\}$. Si l'on note $F^{gr}(V)$ le plus petit sous-espace vectoriel $\Gr_n(V)$-gradué de $\omega_n(X)(V)$
contenant $F(V)$, ce qui précède montre que $C_\infty$ est le plus petit
sous-foncteur de $\omega_n(X)$ tel que $F^{gr}(V)\subset C_\infty(V)$
pour tout $V\in {\rm Ob}\,\E^f$ ; en particulier, $C_\infty\supset F$. 

Soit $Y$ le
plus petit sous-objet de $X$ tel que $Y(V,W)$ contienne les composantes
dans $X(V,W)$ des éléments de $F(V)\subset\omega_n(X)(V)$. La
proposition~\ref{prfig} montre d'une part que $\omega_n(Y)$ est le plus petit sous-$\bar{G}(n)$-comodule de
  $\omega_n(X)$ contenant $F$. Le paragraphe précédent montre d'autre part que
  $C_\infty=\omega_n(Y)$, d'où le troisième point.

Soit $f : C_\infty/F\to\omega_n(T)$ un morphisme de $\F$. Si $X$ est
localement fini, il en est de même pour $Y$, donc la proposition~\ref{fondcr} montre que la composée $g : \omega_n(Y)=C_\infty\twoheadrightarrow
C_\infty/F\xrightarrow{f}\omega_n(T)$ est induite par un morphisme $u : Y\to T$
de $\F_{\Gr,n}$. Comme la composée
$F\hookrightarrow\omega_n(Y)\xrightarrow{g}\omega_n(T)$ est nulle, $F$ est inclus dans
$\omega_n(ker\,u)$. D'après le troisième point, on en déduit $ker\,u=Y$,
puis $f=0$, d'où la dernière assertion.

Par ailleurs, si $X$ est noethérien, alors $Y$ est de type fini, donc $C_\infty=\omega_n(Y)$
est de type fini. Cela montre que la suite $(C_k)_{k>0}$ est
stationnaire, et achève la démonstration.
\end{proof}

On rappelle que la notation  $\F^{lf}_{\Gr,n}$ utilisée dans le
théorème fondamental suivant désigne la
sous-catégorie localisante des
objets localement finis de~$\F_{\Gr,n}$.

\begin{theo}[Théorème de simplicité généralisé]\label{thsg}  Le foncteur 
$$\overline{\omega}_n : \F^{lf}_{\Gr,n}\hookrightarrow\mathcal{F}_{\mathcal{G}r,n}\xrightarrow{\omega_n}\mathcal{F}\twoheadrightarrow\mathcal{F}/\overline{\mathcal{N}il}_{\tilde{\nabla}_n}$$
induit une équivalence entre la catégorie $\F^{lf}_{\Gr,n}$ et une sous--catégorie localisante de
$\overline{\mathcal{N}il}_{\tilde{\nabla}_{n+1}}/\overline{\mathcal{N}il}_{\tilde{\nabla}_n}$. En
particulier, il envoie un foncteur simple de
$\mathcal{F}_{\mathcal{G}r,n}$ sur un objet simple de
$\mathcal{F}/\overline{\mathcal{N}il}_{\tilde{\nabla}_n}$.
\end{theo}

\begin{proof} La proposition~\ref{pr-nno} montre que la restriction à
  $\F_{\Gr,n}^{lf}$ du foncteur $\omega_n$ est bien à valeurs dans
  $\overline{\mathcal{N}il}_{\tilde{\nabla}_{n+1}}$. La proposition~\ref{fondcr} et le théorème~\ref{nabpft} entraînent la pleine fidélité de $\overline{\omega}_n$ et
la stabilité par extensions de son image. Notons $\C_n$ l'image de la
restriction de ce foncteur aux objets {\em finis} de $\F_{\Gr,n}$ :
comme $\overline{\omega}_n$ commute aux colimites, il suffit 
d'établir qu'un sous-objet d'un objet de $\C_n$ est isomorphe à un objet
de $\C_n$. On procède par récurrence sur le degré polynomial. On se donne donc un objet fini $X$
de $\F_{\Gr,n}$ de degré $d\geq 0$ et l'on suppose que l'hypothèse
suivante est satisfaite.

\begin{hyp}[Hypothèse de récurrence] L'image par le foncteur
  $\overline{\omega}_n$ de la sous-catégorie $(\F^f_{\Gr,n})^{d-1}$
  des foncteurs de degré $<d$ de
  $\F^{f}_{\Gr,n}$ est une sous-catégorie épaisse de
  $\F/\overline{\mathcal{N}il}_{\tilde{\nabla}_n}$.
\end{hyp}

\begin{nota} Dans cette section, nous noterons $\A_{n,d}$ cette sous-catégorie épaisse, et $\mathcal{Q}_{n,d}$ la catégorie quotient de
$\F/\overline{\mathcal{N}il}_{\tilde{\nabla}_n}$ par $\A_{n,d}$. Nous noterons également $\mathcal{X}_n$ la sous-catégorie
  pleine des  objets $A$ de $\F$ tels que ${\rm
    hom}\,(A,\omega_n(T))=0$ pour tout objet $T$ de $\F_{\Gr,n}$.
\end{nota}

\begin{lm}\label{lm-prelf5}\begin{enumerate}\item Si $F$ est un objet
    de $\F$ tel que  $\tilde{\nabla}_n (F)$ appartienne à 
$\overline{\mathcal{N}il}_{\tilde{\nabla}_n}$, alors $F$ appartient à $\overline{\mathcal{N}il}_{\tilde{\nabla}_n}$.
\item La sous-catégorie $\mathcal{X}_n$ de $\F$ est stable par le
  foncteur $\tilde{\nabla}_n$.  
\item Un objet de $\mathcal{X}_n$ dont l'image dans
  $\F/\overline{\mathcal{N}il}_{\tilde{\nabla}_n}$ appartient à
  $\A_{n,d}$ est objet de  $\overline{\mathcal{N}il}_{\tilde{\nabla}_n}$.
\item Pour tout entier $k>0$, le morphisme $\pi^k_{n,X}$ induit un
  isomorphisme dans la catégorie $\mathcal{Q}_{n,d}$.
\end{enumerate}
\end{lm}

\begin{proof}[Démonstration du lemme] Soit $F$ est un objet de $\F$ tel que $\tilde{\nabla}_n F$
  appartienne à $\overline{\mathcal{N}il}_{\tilde{\nabla}_n}$. Si $F$ est de type
  fini, $\tilde{\nabla}_n F$ est aussi de type fini (car c'est un quotient
  de $\Delta^n F$), donc $\tilde{\nabla}_n$-nilpotent puisqu'objet de
  $\overline{\mathcal{N}il}_{\tilde{\nabla}_n}$. Par conséquent, $F$
  est $\tilde{\nabla}_n$-nilpotent. Dans le cas général, on montre que
  $F$ est objet de $\overline{\mathcal{N}il}_{\tilde{\nabla}_n}$ en écrivant $F$ comme
  colimite de ses sous-objets de type fini.

Pour le deuxième point, on remarque que la sous-catégorie
$\mathcal{X}_n$ est stable par quotient et préservée par le foncteur $\Delta_{E_n}$, en
raison de  l'isomorphisme ${\rm hom}_\F
(\Delta_{E_n} A,\omega_n(T))\simeq {\rm hom}_\F
(A,\omega_n(T\otimes\iota_n(I_{E_n})))$. 

Le troisième point résulte de la définition de la sous-catégorie
$\A_{n,d}$.

La proposition~\ref{kerpi-cr} montre que le noyau de l'épimorphisme
$\pi_{n,X}$ (corollaire~\ref{crf-nabp}) appartient à
$\A_{n,d}$ (on rappelle que $X$ est de degré $d$). Pour en déduire,
par récurrence sur $k$, la dernière assertion, il suffit de noter
que l'image par le foncteur $\tilde{\nabla}_n$ d'un morphisme $f$ de $\F$
qui induit un isomorphisme dans
$\F/\overline{\mathcal{N}il}_{\tilde{\nabla}_n}$ vérifie encore la
même propriété. Ce résultat s'obtient en notant que les endofoncteurs
{\em exacts} $\Delta$ et $\Delta^n$ de $\F$ préservent
$\overline{\mathcal{N}il}_{\tilde{\nabla}_n}$ (par la dernière
assertion de la proposition \ref{evid-nab}), de sorte que
$ker\,(\Delta f)\simeq\Delta (ker\,f)$, dont $ker\,\tilde{\nabla}_n
(f)$ est un sous-objet, et $coker\,(\Delta^n f)\simeq\Delta^n (coker\,f)$, dont $coker\,\tilde{\nabla}_n
(f)$ est un quotient, sont objets de
$\overline{\mathcal{N}il}_{\tilde{\nabla}_n}$.
\end{proof}

\noindent
{\em Fin de la démonstration du théorème \ref{thsg}.} --- Soit $F$ un sous-objet de $\omega_n(X)$ ; on conserve les notations
de la proposition \ref{cme2}, et l'on se donne $k\in\mathbb{N}^*$ tel
que $C_\infty=C_k$ (qui est donc de la forme $\omega_n(Y)$ pour un sous-objet
$Y$ de $X$). Alors $(\tilde{\nabla}_n)^{k}F$ et $(\tilde{\nabla}_n)^{k}C_\infty$ ont la même
image $C_\infty$ par $\pi_X^k$, qui induit un isomorphisme dans
$\mathcal{Q}_{n,d}$ (dernière assertion du lemme \ref{lm-prelf5}), donc l'inclusion
$(\tilde{\nabla}_n)^{k}F\hookrightarrow(\tilde{\nabla}_n)^{k}C_\infty$ induit un isomorphisme
dans $\mathcal{Q}_{n,d}$. Ainsi, l'image dans
$\F/\overline{\mathcal{N}il}_{\tilde{\nabla}_n}$ de $(\tilde{\nabla}_n)^{k}C_\infty/(\tilde{\nabla}_n)^{k}F$ est objet de
$\A_{n,d}$. Il en est de même pour
$(\tilde{\nabla}_n)^{k}(C_\infty/F)$, qui est un quotient de
$(\tilde{\nabla}_n)^{k}C_\infty/(\tilde{\nabla}_n)^{k}F$
(parce que $\tilde{\nabla}_n$ préserve injections et
surjections). Mais  $(\tilde{\nabla}_n)^{k}(C_\infty/F)$ est aussi un objet de $\mathcal{X}_n$ par la dernière assertion de
la proposition \ref{cme2} et la deuxième assertion du lemme
\ref{lm-prelf5}. Il montre alors que
$(\tilde{\nabla}_n)^{k}(C_\infty/F)$ appartient à
$\overline{\mathcal{N}il}_{\tilde{\nabla}_n}$ (troisième assertion), donc
aussi $C_\infty/F$ (première assertion). Cela achève la démonstration.
\end{proof}

 On démontre de manière similaire
la variante suivante du théorème~\ref{thsg} :

\begin{theo}\label{thsg2} Le foncteur $\omega_n$ induit une équivalence entre la
  sous-catégorie pleine des objets finis de $\F_{\Gr,n}$ et une
  sous-catégorie épaisse de $\N
il_{\tilde{\nabla}_{n+1}}/\N
il_{\tilde{\nabla}_n}$.
\end{theo}

\begin{rem}\label{verbiage2}L'une des conséquences principales de ce théorème
  est le fait que l'image d'un $\bar{G}(n)$-comodule simple dans la
  catégorie quotient  $\F/\N il_{\tilde{\nabla}_n}$ est simple. Pour un
  $\bar{G}(n)$-comodule simple associé à un objet simple
  pseudo-constant de $\F_{\Gr,n}$ (i.e. pour un foncteur de Powell),
  ce résultat est dû à Powell (cf. \cite{GP2}, théorème 6.0.1, et son
  corollaire, la proposition 6.1.1), qui l'a nommé {\em théorème de
    simplicité}. Cela justifie la terminologie employée.

La démonstration du théorème \ref{thsg} repose exactement sur le même
principe que le théorème de simplicité de Powell, à savoir la
considération explicite d'éléments dans les foncteurs, rendue
raisonnable par le calcul aisé du foncteur $\tilde{\nabla}_n$ sur les
projectifs standard, adapté aux
catégories de foncteurs en grassmanniennes par les résultats
préliminaires de la section \ref{par-nablom}.

Le théorème \ref{thsg} (ou~\ref{thsg2}) est un résultat {\em global} sur la structure
de la catégorie $\F$ (il donne des informations sur tous ses
objets de type fini) ; comme le théorème \ref{th-om1}, qui en
constitue le cas particulier $n=1$, il n'utilise pas la théorie des
représentations linéaires (la généralisation à tous les
$\bar{G}(n)$-comodules simples du théorème de simplicité de Powell
fait disparaître les quelques considérations explicites sur les
représentations de $GL_n$ utilisées dans \cite{GP2}).
\end{rem}

\begin{rem} On peut généraliser le théorème~\ref{thsg} à un corps fini
  quelconque, mais il convient de remplacer les foncteurs
  $\tilde{\nabla}_n$ par $\tilde{\nabla}_{(q-1)n}$, où $q$ désigne le
  cardinal du corps considéré. D'autres variantes sont possibles, en
  utilisant également la décomposition scalaire (cf. section~\ref{sct-rqconj},
  pour le cas $n=1$).
\end{rem}

\section{Foncteurs $\tilde{\nabla}_n$-adaptés}\label{sct-ca2}

\begin{conv} Comme dans la section précédente, $n$
désigne un entier strictement positif.
\end{conv}

La notion de foncteur
$\tilde{\nabla}_n$-adapté  que nous introduisons ci-dessous est destinée à faciliter certains
raisonnements de récurrence pour progresser dans l'étude de la
conjecture artinienne extrêmement forte.

\begin{defi}\label{df-ntad} Soit $F$ un objet de $\F$. On
dit que $F$ est un {\em foncteur $\tilde{\nabla}_n$-adapté}
si tout quotient $\tilde{\nabla}_n$-nilpotent de $F$ est oméga-adapté
de hauteur strictement inférieure à~$n$ (cf. §\,\ref{par-omad}).
\end{defi}


Le théorème de simplicité généralisé montre que la conjecture
artinienne extrêmement forte équivaut à dire que tout foncteur de type
fini est $\tilde{\nabla}_i$-adapté pour tout entier $i>0$.

\smallskip

La définition~\ref{df-ntad} est très difficile à vérifier si l'on
ignore si un quotient d'un foncteur oméga-adapté
de hauteur strictement inférieure à $n$ est encore oméga-adapté
de hauteur strictement inférieure à $n$ : elle est maniable lorsque l'hypothèse suivante
est satisfaite (le but étant de dé\-montrer l'épaisseur de
$\F^{\omega-ad(n)}$).

\begin{hyp}\label{htca} La sous-catégorie $\F^{\omega-ad(i)}$ de
  $\F$ est épaisse pour $i<n$.
\end{hyp}

\begin{pr}\label{pr-fnta} Soient $A\in {\rm
    Ob}\,\F$ et $F$ un objet {\bf fini} de $\F$. On suppose que
  l'hypothèse \ref{htca} est satisfaite.
\begin{enumerate}\item Si $A$ est un foncteur
  $\tilde{\nabla}_n$-adapté, il en est de même pour tous ses quotients.
\item  Si $A$ est un foncteur $\tilde{\nabla}_n$-adapté, alors $A\otimes F$
  est  $\tilde{\nabla}_n$-adapté.
\item Si $A$ est  $\tilde{\nabla}_n$-adapté, alors $(A : F)$ est  $\tilde{\nabla}_n$-adapté.
\end{enumerate}
\end{pr}

\begin{proof} Le premier point est formel. 

Pour le second, considérons un épimorphisme $ f :A\otimes
F\twoheadrightarrow Q$, où $Q$ est $\tilde{\nabla}_n$-nilpotent. Soit
$g : A\to\mathbf{Hom}_\F (F,Q)$ le morphisme adjoint à $f$. Comme $F$
est quotient d'une somme directe finie de projectifs standard,
$\mathbf{Hom}_\F (F,Q)$ est un sous-foncteur d'une somme directe finie
de $\Delta_V Q$, il est donc $\tilde{\nabla}_n$-nilpotent par la
proposition \ref{evid-nab}.\,\ref{p4n}. Par conséquent, $im\,g$ est un
quotient $\tilde{\nabla}_n$-nilpotent de $A$, c'est donc, par
hypothèse, un foncteur oméga-adapté
de hauteur strictement inférieure à $n$. Il en est de même pour
$im\,g\otimes F$, puisque $F$ est fini (cf. \cite{art2}, §\,$12.1$).

Le diagramme commutatif
$$\xymatrix{ & A\otimes
  F\ar@{>>}[d]^-f\ar[dl]_-{g\otimes F} \\
\mathbf{Hom}_\F (F,Q)\otimes F\ar[r] & Q
}$$
(dont la flèche horizontale est la coünité de
l'adjonction) montre que $Q$
est un quotient de $im\,g\otimes F$ ; $Q$ est donc objet de
$\F^{\omega-ad(n-1)}$ grâce à ce qui précède et à l'hypothèse~\ref{htca}.

Le troisième point s'établit par un argument d'adjonction analogue, grâce aux deux remarques suivantes.
\begin{enumerate}\item Le foncteur $\cdot\otimes F$ préserve $\N
  il_{\tilde{\nabla}_n}$, par la proposition \ref{nab-tens}.
\item Le foncteur de division par $F$ préserve
  $\F^{\omega-ad (n-1)}$. En effet, comme $F$ est de co-type fini, ce
  foncteur est un quotient d'une somme directe finie de foncteurs
de  décalage, et la sous-catégorie $\F^{\omega-ad (n-1)}$ de $\F$ est
  stable par les foncteurs de décalage (cf. \cite{art2}, §\,$12.1$).
\end{enumerate} 
\end{proof}

Avant d'appliquer cette propriété à la proposition~\ref{reduca}, nous mentionnons un lemme élémentaire qui se déduit
des résultats de \cite{K2},~§\,$4$.

\begin{lm}\label{tec-rg} Soit $M$ un
  $GL_n$-module fini. Il existe un foncteur {\bf fini} $F$ de $\F$ tel
  que $F(E_n)$ est isomorphe à $M$ comme $GL_n$-module.
\end{lm}

La proposition suivante permet de réduire les vérifications
nécessaires pour démontrer que $\F^{\omega-ad(n)}$ est épaisse.

\begin{pr}\label{reduca} Soit $\lambda$ une partition
  régulière telle que $\lambda_1=n$. On suppose l'hypothèse~\ref{htca} vérifiée. Les assertions suivantes sont
  équivalentes.
\begin{enumerate}\item\label{ire1} La sous-catégorie $\F^{\omega-ad(n)}$ de $\F$
  est épaisse.
\item\label{ire2} Pour tout objet fini $X$ de $\F_{\Gr,n}$, le foncteur
  $\omega_n(X)$ est $\tilde{\nabla}_n$-adapté.
\item\label{ire3} Pour tout objet simple $S$ de $\F_{\Gr,n}$, le foncteur
  $\omega_n(S)$ est $\tilde{\nabla}_n$-adapté.
\item\label{ire4} Le foncteur de Powell $Q_\lambda$ est $\tilde{\nabla}_n$-adapté.
\end{enumerate}

Lorsqu'elles sont vérifiées, pour tout foncteur fini (resp. simple)
$X$ de $\F_{\Gr,k}$, où $k\leq n$, le foncteur $\omega_k(X)$ est
noethérien (resp. simple noethérien) de type $k$. 
\end{pr}

\begin{proof} Supposons l'assertion~\ref{ire4} vérifiée. On
  commence par montrer que pour tout $GL_n$-module simple $S$, le
  foncteur de Powell $\omega_n\rho_n(S)$ est
  $\tilde{\nabla}_n$-adapté. Pour cela, on se donne, conformément au
  lemme \ref{tec-rg}, un foncteur fini $F$ de $\F$ tel que le
  $GL_n$-module $F(E_n)$ est isomorphe à $\FF[GL_n]$. Le foncteur
  $\omega_n(\rho_n(R_\lambda) : \iota_n(F))$ est $\tilde{\nabla}_n$-adapté,
  car c'est un quotient, par la proposition~$9.2$ de \cite{art2}, de
  $\omega_{\leq n}(\mathcal{P}_{n,\leq n}\,\rho_n(R_\lambda) : \iota_{\leq
    n}(F))\simeq (\omega_n\rho_n(R_\lambda) : F)=(Q_\lambda : F)$ (cet isomorphisme
  venant de la proposition~$9.8$ de \cite{art2}), de sorte que la proposition~\ref{pr-fnta} prouve ce premier point.

Par ailleurs, $(\rho_n(R_\lambda) : \iota_n(F))\simeq\rho_n(R_\lambda :
\FF[GL_n])$ par la proposition \ref{divrho}, et $(R_\lambda :
\FF[GL_n])\simeq R_\lambda\otimes\FF[GL_n]\simeq\FF[GL_n]^{\oplus i}$,
où $i=\dim_\FF R_\lambda\in\mathbb{N}^*$. Par conséquent, tout $GL_n$-module simple $S$ est quotient
de  $(R_\lambda : \FF[GL_n])$, donc le quotient $\omega_n\rho_n(S)$ de
$\omega_n(\rho_n(R_\lambda) : \iota_n(F))$ est
$\tilde{\nabla}_n$-adapté.

Si $X$ est un objet simple de $\F_{\Gr,n}$, il existe un $GL_n$-module simple
$S$ et un objet simple $F$ de $\F$ tels que
$X\simeq\kappa_n(F)\otimes\rho_n(S)$ (proposition \ref{prcqd}), donc
$\omega_n(X)$ est quotient de
$\omega_n(\iota_n(F)\otimes\rho_n(S))\simeq \omega_n\rho_n(S)\otimes
F$. La proposition \ref{pr-fnta} montre à nouveau que ce foncteur est
$\tilde{\nabla}_n$-adapté. 

On a ainsi démontré que \ref{ire4} implique~\ref{ire3}.

\smallskip

Si l'assertion \ref{ire3} est vérifiée, le théorème~\ref{thsg}
prouve que tout sous-quotient de $\omega_n(S)$, où $S$ est un objet
simple de $\F_{\Gr,n}$, est objet de $\F^{\omega-ad(n)}$. Comme la
sous-catégorie pleine $\A$ des objets de $\F$ dont tous les sous-quotients
sont dans $\F^{\omega-ad(n)}$ vérifie l'hypothèse que pour toute suite exacte courte $0\to A\to B\to
  C\to 0$ de $\F$, si deux des objets $A, B, C$ appartiennent à $\A$, il en
  est de même du troisième, on en déduit l'assertion~\ref{ire1}.

Il est clair que~\ref{ire2} entraîne~\ref{ire4}.

Si l'assertion~\ref{ire1} est vérifiée, le foncteur $\omega_n$
induit une équivalence entre les catégories $\F^f_{\Gr,n}$ et
$\F^{\omega - ad(n)}/\F^{\omega - ad(n-1)}$, donc l'assertion~\ref{ire2} est satisfaite. 

Ainsi, les assertions de l'énoncé sont équivalentes.

La fin de la proposition découle de la proposition~\ref{rec-omad}.
\end{proof}

\section{Structure de $P^{\otimes 2}\otimes F$ pour un
  foncteur fini $F$}\label{sct-concr}

Dans l'article \cite{GP5} (proposition~$7.4$) est établi le résultat
suivant\,\footnote{On rappelle que l'on travaille avec des conventions
duales de celles de Powell.}.

\begin{pr}[Powell]\label{prngp} Soit $F$  un foncteur de type fini de $\F$
  tel que $\tilde{\nabla}_2(F)=0$ et que ${\rm
    hom}_\F(F,\bar{P})=0$. Alors $F$ est un foncteur fini.
\end{pr}

En étudiant une filtration explicite du foncteur $\bar{G}(2)$
construite à partir de la filtration co-polynomiale de $\bar{P}$,
\cite{GP5} déduit de la proposition précédente le corollaire suivant
(c'est sa proposition~$7.5$).

\begin{cor}[Powell]\label{crngp} Le foncteur $\bar{G}(2)$ est $\tilde{\nabla}_2$-adapté.
\end{cor}

Nous aurons également besoin de la conséquence facile suivante de
la proposition~\ref{prngp} (qui est implicite dans \cite{GP5}).

\begin{cor}\label{cr-dgp} Si $F$ est un foncteur de type fini de $\F$
  tel que $\tilde{\nabla}_2(F)=0$, alors $F$ est oméga-adapté de
  hauteur au plus~$1$. Plus précisément, $F$ est isomorphe dans
  $\F/\F^f$ à une somme directe finie de copies de~$\bar{P}$.
\end{cor}

\begin{proof} Comme $F$ est de type fini, l'ensemble ${\rm
    hom}_\F(F,\bar{P})$ est fini. Soient $N$ et $C$ le noyau et le
  conoyau, respectivement, du morphisme canonique $F\to\bar{P}^{{\rm
    hom}_\F(F,\bar{P})}$. Alors $C$ est fini (car ${\rm
  hom}_\F(C,\bar{P})=0$) ; comme $\bar{P}$ est $\F^f$-parfait, on en
déduit ${\rm
  hom}_\F(N,\bar{P})=0$. Mais $\tilde{\nabla}_2(N)=0$, car $N$ est un
sous-foncteur de $F$, donc $N$ est fini par la
proposition~\ref{prngp}, d'où le corollaire.
\end{proof}

\begin{rem} La démonstration de la proposition~\ref{prngp} repose sur
  des considérations explicites issues de la théorie des
  représentations --- essentiellement, Powell s'appuie sur le fait
  qu'il existe peu d'extensions non triviales entre les puissances
  extérieures (cf. \cite{Franjou}), qui sont les seuls foncteurs simples de $\F$ annihilés
  par $\tilde{\nabla}_2$. Malheureusement, la complexité des problèmes
  posés par la généralisation de cette approche aux foncteurs
  $\tilde{\nabla}_n$ croît très rapidement avec~$n$.

Du reste, la seule obstruction sérieuse à la
généralisation des résultats présentés dans cette section se trouve
concentrée à la proposition~\ref{prngp}.
\end{rem}

La satisfaction de l'hypothèse \ref{htca} de la section précédente pour $n=2$
(corollaire~\ref{th1-oad}) et le corollaire~\ref{crngp} permettent de
déduire de la proposition~\ref{reduca} le théorème suivant, qui
constitue le résultat~\guillemotleft~concret~\guillemotright~le plus important de cet article.

\begin{theo}\label{avca2}\begin{enumerate}\item La sous-catégorie
    $\F^{\omega-ad(2)}$ de $\F$ est épaisse.
\item Tout $\bar{G}(2)$-comodule fidèle fini (resp. simple) est noethérien
  (resp. simple noethérien) de type $2$.
\item En particulier, pour tout objet fini $F$ de $\F$, le foncteur
  $P^{\otimes 2}\otimes F$ est noethérien.
\end{enumerate}
\end{theo}

Ce théorème semble le
meilleur résultat actuellement connu concernant la conjecture
artinienne. Dans le cas où $F$ est constant, il est dû à Powell
(cf. \cite{GP5}) ; nous avons traité le cas où $F$ est une puissance
extérieure  par d'autres méthodes dans \cite{art1}. 
\smallskip

Nous pouvons maintenant préciser le corollaire~\ref{cr-dgp} en montrant que
  la conjecture~\ref{canabla} est vérifiée pour $n=2$.

\begin{pr}\label{onnapbse} Tout foncteur $\tilde{\nabla}_2$-nilpotent
  et de type fini $F$ est oméga-adapté de
  hauteur au plus~$1$.
\end{pr}

\begin{proof} On établit la proposition par récurrence sur l'indice de
  $\tilde{\nabla}_2$-nilpotence $i$ de $F$. On suppose donc la
  propriété vérifiée pour tous les foncteurs de type fini annihilés
  par $(\tilde{\nabla}_2)^{i-1}$.

Il existe
une suite exacte
$$\bar{P}\otimes\tilde{\nabla}_2(F)\to F\to Q\to 0$$
où $\tilde{\nabla}_2(Q)=0$ (cf.~\cite{GP4}, §\,$5$). Comme $Q$ est de type
fini, on en déduit que $Q$ est oméga-adapté de
  hauteur au plus~$1$ par le corollaire~\ref{cr-dgp}. Le foncteur de
  type fini $\tilde{\nabla}_2(F)$ étant annulé par
  $(\tilde{\nabla}_2)^{i-1}$, l'hypothèse de récurrence montre qu'il
  est également oméga-adapté de
  hauteur au plus~$1$, donc $\bar{P}\otimes\tilde{\nabla}_2(F)$, puis $F$, sont oméga-adaptés de
  hauteur au plus~$2$ (cf. \cite{art2}, §\,$12.1$ et le théorème~\ref{avca2}). En conséquence, comme le foncteur $\omega_2$
  induit une équivalence entre les catégories $\F^f_{\Gr,2}$ et
  $\F^{\omega - ad(2)}/\F^{\omega - ad(1)}$, il existe une suite exacte
$$0\to A\to F\to\omega_2(X)\to B\to 0$$
avec $X\in {\rm Ob}\,\F^f_{\Gr,2}$ et $A, B\in {\rm Ob}\,\F^{\omega -
  ad(1)}$. Vu que la sous-catégorie $\N il_{\tilde{\nabla}_2}$ de $\F$ est
épaisse (proposition~\ref{nabnil-ep}) et que $A$, $B$ et $F$ sont
$\tilde{\nabla}_2$-nilpotents, cela entraîne que $\omega_2(X)$ est
$\tilde{\nabla}_2$-nilpotent, donc que $X=0$ par la
proposition~\ref{nabpft}. 

Par conséquent, $F$ est oméga-adapté de hauteur au plus~$1$, ce qu'il
fallait démontrer.
\end{proof}


\subsection*{Remerciements} L'auteur exprime sa chaleureuse gratitude
envers Geoffrey Powell pour l'attention qu'il a portée à cet
article. Il remercie également Lionel Schwartz pour ses conseils.

\nocite{*}
\bibliographystyle{smfalpha}
\bibliography{bibart3}

\providecommand{\bysame}{\leavevmode ---\ }
\providecommand{\og}{``}
\providecommand{\fg}{''}
\providecommand{\smfandname}{\&}
\providecommand{\smfedsname}{\'eds.}
\providecommand{\smfedname}{\'ed.}
\providecommand{\smfmastersthesisname}{M\'emoire}
\providecommand{\smfphdthesisname}{Th\`ese}
\begin{thebibliography}{FFPS03}

\bibitem[CR87]{CR2}
{\scshape C.~W. Curtis {\normalfont \smfandname} I.~Reiner} -- \emph{Methods of
  representation theory. {V}ol. {II}}, Pure and Applied Mathematics (New York),
  John Wiley \& Sons Inc., New York, 1987, With applications to finite groups
  and orders, A Wiley-Interscience Publication.

\bibitem[CR90]{CR}
\bysame , \emph{Methods of representation theory. {V}ol. {I}}, Wiley Classics
  Library, John Wiley \& Sons Inc., New York, 1990, With applications to finite
  groups and orders, Reprint of the 1981 original, A Wiley-Interscience
  Publication.

\bibitem[Dja]{these}
{\scshape A.~Djament} -- {\og Repr\'esentations g\'en\'eriques des groupes
  lin\'eaires : cat\'egories de foncteurs en grassmanniennes, avec applications
  à la conjecture artinienne\fg}, \smfphdthesisname, Universit\'e Paris 13, en
  pr\'eparation.

\bibitem[Dja06a]{art2}
\bysame , {\og Cat\'egories de foncteurs en grassmanniennes\fg},
  arXiv:math.AT/0610598, 2006.

\bibitem[Dja06b]{art1}
\bysame , {\og Foncteurs de division et structure de {$I^{\otimes
  2}\otimes\Lambda^n$} dans la cat\'egorie {$\mathcal{F}$}\fg},
  arXiv:math.RT/0607595, 2006.

\bibitem[FFPS03]{FFPS}
{\scshape V.~Franjou, E.~M. Friedlander, T.~Pirashvili {\normalfont
  \smfandname} L.~Schwartz} -- \emph{Rational representations, the {S}teenrod
  algebra and functor homology}, Panoramas et Synth\`eses [Panoramas and
  Syntheses], vol.~16, Soci\'et\'e Math\'ematique de France, Paris, 2003.

\bibitem[Fra96]{Franjou}
{\scshape V.~Franjou} -- {\og Extensions entre puissances ext\'erieures et
  entre puissances sym\'etriques\fg}, \emph{J. Algebra} \textbf{179} (1996),
  no.~2, p.~501--522.

\bibitem[Gab62]{Gab}
{\scshape P.~Gabriel} -- {\og Des cat\'egories ab\'eliennes\fg}, \emph{Bull.
  Soc. Math. France} \textbf{90} (1962), p.~323--448.

\bibitem[Hum87]{Hum}
{\scshape J.~E. Humphreys} -- {\og The {S}teinberg representation\fg},
  \emph{Bull. Amer. Math. Soc. (N.S.)} \textbf{16} (1987), no.~2, p.~247--263.

\bibitem[Jam78]{James}
{\scshape G.~D. James} -- \emph{The representation theory of the symmetric
  groups}, Lecture Notes in Mathematics, vol. 682, Springer, Berlin, 1978.

\bibitem[Jan03]{Jan}
{\scshape J.~C. Jantzen} -- \emph{Representations of algebraic groups}, second
  \smfedname, Mathematical Surveys and Monographs, vol. 107, American
  Mathematical Society, Providence, RI, 2003.

\bibitem[Kuh94a]{K1}
{\scshape N.~J. Kuhn} -- {\og Generic representations of the finite general
  linear groups and the {S}teenrod algebra. {I}\fg}, \emph{Amer. J. Math.}
  \textbf{116} (1994), no.~2, p.~327--360.

\bibitem[Kuh94b]{K2}
\bysame , {\og Generic representations of the finite general linear groups and
  the {S}teenrod algebra. {II}\fg}, \emph{$K$-Theory} \textbf{8} (1994), no.~4,
  p.~395--428.

\bibitem[Mit86]{Mitch}
{\scshape S.~Mitchell} -- {\og On the {S}teinberg module, representations of
  the symmetric groups, and the {S}teenrod algebra\fg}, \emph{J. Pure Appl.
  Algebra} \textbf{39} (1986), no.~3, p.~275--281.

\bibitem[Pir97]{Piriou}
{\scshape L.~Piriou} -- {\og Sous-objets de {$\overline I\otimes \Lambda\sp n$}
  dans la cat\'egorie des foncteurs entre {$\bold F\sb 2$}-espaces
  vectoriels\fg}, \emph{J. Algebra} \textbf{194} (1997), no.~1, p.~53--78.

\bibitem[Pop73]{Pop}
{\scshape N.~Popescu} -- \emph{Abelian categories with applications to rings
  and modules}, Academic Press, London, 1973, London Mathematical Society
  Monographs, No. 3.

\bibitem[Pow98a]{GP5}
{\scshape G.~M.~L. Powell} -- {\og The {A}rtinian conjecture for {$I\sp
  {\otimes2}$}\fg}, \emph{J. Pure Appl. Algebra} \textbf{128} (1998), no.~3,
  p.~291--310, With an appendix by Lionel Schwartz.

\bibitem[Pow98b]{GP4}
\bysame , {\og Polynomial filtrations and {L}annes' {$T$}-functor\fg},
  \emph{$K$-Theory} \textbf{13} (1998), no.~3, p.~279--304.

\bibitem[Pow98c]{GP2}
\bysame , {\og The structure of indecomposable injectives in generic
  representation theory\fg}, \emph{Trans. Amer. Math. Soc.} \textbf{350}
  (1998), no.~10, p.~4167--4193.

\bibitem[Pow00a]{GP1}
\bysame , {\og On {A}rtinian objects in the category of functors between
  {$\bold F\sb 2$}-vector spaces\fg}, in \emph{Infinite length modules
  (Bielefeld, 1998)}, Trends Math., Birkh\"auser, Basel, 2000, p.~213--228.

\bibitem[Pow00b]{GP3}
\bysame , {\og The structure of the tensor product of {$\bold F\sb 2[-]$} with
  a finite functor between {$\bold F\sb 2$}-vector spaces\fg}, \emph{Ann. Inst.
  Fourier (Grenoble)} \textbf{50} (2000), no.~3, p.~781--805.

\bibitem[Pow01]{GP-dim}
\bysame , {\og The tensor product theorem for {$\widetilde\nabla$}-nilpotence
  and the dimension of unstable modules\fg}, \emph{Math. Proc. Cambridge
  Philos. Soc.} \textbf{130} (2001), no.~3, p.~427--439.

\bibitem[PS98]{PS1}
{\scshape L.~Piriou {\normalfont \smfandname} L.~Schwartz} -- {\og Extensions
  de foncteurs simples\fg}, \emph{$K$-Theory} \textbf{15} (1998), no.~3,
  p.~269--291.

\end{thebibliography}

\end{document}